\documentclass[journal,twocolumn]{IEEEtran}

\ifCLASSINFOpdf
\else
\fi

\usepackage{url,amsmath,amsfonts,amssymb,bm}
\usepackage{pifont}
\usepackage{enumitem}
\usepackage{graphicx}
\usepackage{algorithm}
\usepackage{algorithmicx, algpseudocode }
\usepackage{multirow}
\usepackage{color}
\usepackage{mathtools}

\usepackage{caption}
\usepackage{subcaption}

\usepackage{amsthm}

\usepackage[colorlinks=true, linkcolor=blue, citecolor=blue, urlcolor=black]{hyperref}
\hypersetup{pdfcreator  = {LaTeX with hyperref package},
	pdfproducer = {dvips + ps2pdf} }
\newtheorem{definition}{Definition}

\newtheorem{theorem}{Theorem}
\newtheorem{lemma}{Lemma}

\newtheorem{proposition}{Proposition}
\newtheorem{assumption}{Assumption}

\newcommand{\E}{{\mathbb E}}

\newcommand{\M}{\mathcal{E}}

  \DeclareMathOperator*{\erank}{r}
  \DeclareMathOperator*{\rank}{{rank}}
    \DeclareMathOperator*{\srank}{r}

  \DeclareMathOperator*{\tr}{tr}
    \DeclareMathOperator*{\diag}{diag}

\DeclareMathOperator*{\var}{var}
\DeclareMathOperator*{\Cov}{Cov}
\DeclareMathOperator*{\Prob}{\mathbb{P}}

\DeclareMathOperator*{\Sy}{\mathbf{S}_{\mathbf{y}}}
\DeclareMathOperator*{\Cy}{\mathbf{C}_{\mathbf{y}}}

\DeclareMathOperator*{\Ry}{\mathbf{R}_{\mathbf{y}}}
\def\NEW#1{\textcolor[rgb]{0.00,0.0,0.0}{#1}}%

\begin{document}
%
%
\title{ Covariance Matrix Estimation with Non Uniform and Data Dependent Missing Observations}
%
%

\author{ Eduardo~Pavez ~\IEEEmembership{Member,~IEEE,}
        and~Antonio~Ortega~\IEEEmembership{Fellow,~IEEE} 
  %
\thanks{Authors are with the Department of Electrical and Computer Engineering,
University of Southern California, Los Angeles,
CA, 90089 USA.  This work was funded in part by NSF  under grants CCF-1410009 and CCF-2009032.
(Author's e-mail: pavezcar@usc.edu, ortega@sipi.usc.edu)}
}

%
%

%


\maketitle
%
\begin{abstract}
In this paper we study covariance estimation with missing data. We consider    missing data mechanisms that can be independent of the data,  or have a time varying   dependency. Additionally,  observed  variables may have arbitrary (non uniform) and dependent observation probabilities.  For each mechanism, we construct an  unbiased estimator and  obtain bounds for the expected value of their estimation error  in operator norm. Our  bounds are equivalent, up to constant and logarithmic factors, to state of the art bounds for complete and uniform missing observations. Furthermore,  for the more general non uniform and dependent cases, the proposed bounds are new or improve upon previous results. Our error estimates depend on quantities we call {\em scaled effective rank}, which generalize the effective rank to account for missing observations. All the estimators studied in this work  have the same asymptotic convergence rate (up to logarithmic factors). 
\end{abstract}
%
%
%
%
%
%
%
%
\IEEEpeerreviewmaketitle
\section{Introduction}
\IEEEPARstart{T}{he} covariance matrix  is an essential component of many algorithms in  machine learning \cite{zhang2010learning}, medical data analysis \cite{schafer2005shrinkage, bullmore_complex_2009},  finance \cite{bai2011estimating},  and signal processing \cite{stoica2011spice,dabov2009bm3d,pyatykh2013image, chen2010shrinkage,egilmez2017graph}. 
Estimation of {covariance} matrices involves the design and analysis of statistical procedures for recovering the  covariance matrix from data samples.
 %
 The prevalence and inevitability of missing data, along with the importance of covariance matrices in many applications and algorithms, makes the study of   covariance estimation  with missing observations  of great importance. 
 %
  
 When there are complete observations, the sample covariance matrix is the  most commonly used estimator. Its performance  is well understood under various settings (e.g., noise, high dimensions). A favored  criteria to evaluate an estimator's performance is  the estimation error measured in operator norm (largest singular value). This norm is preferred over other matrix norms (Frobenius, Schatten, etc.) because    bounds for   other norms   imply operator norm bounds, while the converse is not true in general. The best  estimation error bounds for the sample covariance matrix  were recently  obtained under sub-Gaussian \cite{lounici2014high,bunea2015sample} and Gaussian \cite{koltchinskii2017concentration} assumptions. Essentially,  these results  indicate that the number of samples required for accurate covariance estimation  is proportional to the  {\em effective rank} of the population covariance matrix,  a parameter that is upper bounded by the  rank, but that can be much smaller. 
 
 Even though missing data is pervasive, and covariance matrices are utilized by a plethora of algorithms,  studies on  covariance estimation with missing data are few and of  a narrower scope than those considering the  complete data case.  Perhaps because there are many available estimators, and each of them is designed and analyzed  for a different type of missing data mechanism. In addition,  mathematical expressions involving these estimators are much more intricate than those for complete observations. The  works of \cite{lounici2014high,cai2016minimax} stand out since they offer error analyses and  bounds that match similar ones obtained for complete data. 
  They  studied convergence to the population covariance matrix (consistency) of different  estimators by deriving finite sample  error bounds. These works assume that the observations are   independent and identically distributed (i.i.d.)  copies of a  sub-Gaussian vector, while the population covariance matrix may be structured or unstructured. Sparse and   bandable structures were considered in \cite{cai2016minimax}, while \cite{lounici2014high} considered both structured (low rank) and unstructured covariance matrices. In addition, both papers assume the missing data mechanism is independent of the observations, however the specific assumptions and estimators are considerably  different. More recently, \cite{park2019non}  considered a non uniform and dependent  missing data pattern. This work  however, when simplified to the case studied by \cite{lounici2014high} (uniform independent observations),  returns  convergence rates that are sub-optimal with respect to the ambient dimension and the rate of missing entries. 

In this work we do not make  any structural assumptions (such as low rank or sparsity) on the population covariance matrix. We consider  various missing data mechanisms and estimators, including {a class of} time varying dependent observations.  Each of these mechanisms has a corresponding unbiased covariance estimator that generalizes the sample covariance matrix. 
For each  case  we provide bounds for the expected value of  estimation error  in operator norm. Our principal contributions are:
\begin{itemize}
	\item We study existing  covariance estimators under  the  {\em  missing completely at random (MCAR)} missing data mechanism \cite{little2014statistical}, which is data  independent  (see Section \ref{sec:miss_estimator}). We consider an instance of this  mechanism in which variables are observed according to {\em dependent} Bernoulli random variables with  {\em non uniform} probabilities (the same missing data mechanism and estimator considered in \cite{park2019non}). We are the first to provide rate optimal error bounds (Theorem \ref{th_expected_error}) for this setting. The case of equal observation probabilities (uniform case) was considered in \cite{lounici2014high}, while the non uniform  case with dependent observation variables was studied in \cite{park2019non}. 
	\item \NEW{We consider another instance of the MCAR mechanism, in which the missing data distribution is unknown, therefore the observation probabilities have to be estimated from data. We propose an unbiased  estimator that is closely related to the one used for estimation of structured covariance matrices  \cite{cai2016minimax,kolar2012consistent}. This  work is the first in obtaining estimation error bounds in operator norm in  the unstructured case (Theorem \ref{th_expected_error_MAR}).  }
	\item We propose the conditionally MCAR (CMCAR) missing data mechanism, in which the missing data distribution evolves over time and depends on previous observations. We propose an unbiased estimator and characterize its estimation error  (Theorem \ref{th_expected_error_martin}). 
	\item  Under certain conditions,  we show that the   bounds from Theorems \ref{th_expected_error_MAR} and \ref{th_expected_error_martin},  can be simplified to resemble the bound from Theorem \ref{th_expected_error}.  This is possible because they all depend on  quantities we call {\em scaled effective rank}, which generalize the effective rank to account for incomplete observations. In addition, after further  simplifications, our bounds match or improve previous state of the art results with complete and missing observations.
\end{itemize}
This paper is organized as follows. In Section \ref{sec:miss_related_cov} we review the literature on covariance estimation with missing data. Section \ref{sec:miss_estimator} introduces the missing data mechanisms and estimators. Error bounds  are presented and discussed in Section \ref{sec:miss_missing}, while numerical experiments are contained in Section \ref{sec:miss_exp_cov}. Proofs and Conclusion appear in Sections \ref{sec:miss_proofs} and \ref{sec:miss_conc} respectively.
\section{Related work}
 \label{sec:miss_related_cov}
\subsection{Applications of covariance matrices}
 Missing data problems  motivating this work are listed next. 
 \paragraph{``Plug-in" estimators}
 Numerous algorithms  use   the sample covariance matrix as input. This  ``plug-in" principle can be used to  extend these techniques to the missing data case  by using  estimators from this work that account for missing data. Examples that have followed this strategy include sparse sub-space clustering  \cite{yang2015sparse},  multi-task learning \cite{hunt2018multi}, sub-space learning \cite{gonen2016subspace},  classification \cite{chi2013nearest}, image compression \cite{chao2017thesis}, regression \cite{loh2012high}, principal component analysis  \cite{jurczak2017spectral,lounici2013sparse},    and inverse covariance estimation \cite{stadler2012missing, kolar2012consistent}.
 \paragraph{Applications with data acquisition  costs} For some applications, obtaining complete observations might not be possible due to resource constraints.  For example, in clinical studies, measuring certain features (e.g., acquiring samples or running tests)  has monetary costs  \cite{veeramachaneni2005active,lim2005hybrid}, while in sensor networks   obtaining measurements incurs in  energy consumption  \cite{asif2016matrix, deshpande2004model,gershenfeld2010intelligent}.  For this reason,  data may be acquired  based on the cost of acquiring each sample, as well as information from previous observations. In these examples a  covariance matrix may be required for statistical analysis, but standard covariance estimators are not well suited since observations are incomplete, possibly non i.i.d., and the missing data mechanism may evolve over time.
 \paragraph{Statistical estimation and inference with partial observations} %
 A very high ambient dimension  can be problematic for data storage,  transmission,  or processing. But, since most natural signals are very structured, a few partial observations are often sufficient for  performing simple statistical tasks like estimation and detection \cite{davenport_signal_2010}. Covariance estimation  from compressed measurements  has been  studied for structured  \cite{dasarathy2015sketching, chen2015exact}, and unstructured covariance matrices \cite{chen2017toward, azizyan2018extreme,pourkamali2017preconditioned,pourkamali2016estimation, anaraki2014memory}. Additionally,  active learning approaches, have been shown to benefit sequential hypothesis testing \cite{naghshvar2013active},  Gaussian graphical model selection  \cite{dasarathy2016active,scarlett2017lower}, and  covariance estimation \cite{pavez2018active}.
 %
%
\subsection{Covariance estimation with complete observations }
The complete data is denoted by the $n \times N$  matrix $\mathbf{X}$, where the columns of $\mathbf{  X }$ are assumed to be zero mean,  independent and  identically distributed (i.i.d.) realizations of a vector $\mathbf{x}$ with population covariance matrix $\mathbf{  \Sigma }$. The starting point of most covariance estimation studies is the sample covariance matrix
given by
\begin{equation}\label{eq_sample_cov_first}
\mathbf{S} = \frac{1}{N} \mathbf{ X X^{\top}},
\end{equation}
which is an unbiased estimator for the population covariance matrix $\mathbf{\Sigma}$, that is, $\E[\mathbf{S}] = \mathbf{\Sigma}$. 
 Some recent papers have shown that the  estimation  error in operator norm of the sample covariance matrix is characterized by   the {\em effective rank} of  $\mathbf{\Sigma}$ \cite{lounici2014high,bunea2015sample,koltchinskii2017concentration}.  
These results imply the expectation bound 
\begin{align}\label{eq_error_cov_complete}
\frac{\E\left[ \Vert  \mathbf{\Sigma -S} \Vert  \right]}{\Vert \mathbf{\Sigma} \Vert} \leq \mathcal{O}\left( \sqrt{\frac{\erank(\mathbf{\Sigma})}{N}} \vee \frac{\erank(\mathbf{\Sigma})}{N}\right),
\end{align}
where $a\vee b$ stands for $\max(a,b)$. The $\mathcal{O}$ notation hides logarithmic dependencies on the dimension $n$ and number of samples $N$. The effective rank $\erank(\mathbf{\Sigma})$ (see Definition \ref{def_erank})    is upper bounded by the actual   rank, but can  be much smaller, thus offering a more nuanced measure of dimensionality.  

Inequality (\ref{eq_error_cov_complete}) implies that the sample complexity of (\ref{eq_sample_cov_first})  is $N=\mathcal{O}(\erank(\mathbf{\Sigma})/ \epsilon^2)$, that is, the number of i.i.d. samples required for the error in (\ref{eq_error_cov_complete}) to be below $\epsilon$ is proportional to the effective rank. 
Results from  \cite{bunea2015sample,lounici2014high} assume a sub-Gaussian distribution, while   \cite{koltchinskii2017concentration} provides dimension free bounds (without logarithmic dependence on $n$) for the Gaussian case. Our bounds are also based on the effective rank, so they can be compared with  results  that assume complete observations.
\subsection{Missing data mechanisms}
%

Construction of unbiased estimators requires knowledge of the {\em missing data mechanism}, i.e., the process
by which  data is lost. Recent literature of covariance estimation with missing  observations considers  mechanisms that  are      \emph{data independent}. This type of mechanism is  called missing completely at random (MCAR) \cite{little2014statistical}. We also   introduce a new one called conditionally MCAR (CMCAR).

The MCAR mechanism has been used to model missing features in machine learning applications \cite{marlin2008missing}, or data loss in  sensor networks, where  transmission errors (e.g., packet loss), 
or sensor failures induce missing data. This mechanism has also been used to  model incomplete observations in clinical and  longitudinal studies \cite{laird1988missing},  in which data is lost due to   subject dropout. 
When the missing data is MCAR, we study two scenarios, {\em known} and {\em unknown} missing data distribution. In the first case, coordinates are observed according to dependent Bernoulli  random variables that  may have different probabilities. A value $0/1$ indicates a  missing/observed variable. Covariance estimators designed for this setting appear in \cite{lounici2014high, loh2012high, park2019non}. In \cite{lounici2014high}, the missing data mechanism is modeled by  independent Bernoulli variables with uniform distribution (all coordinates are observed with the same probability), while \cite{park2019non} addresses the dependent non uniform case. In  \cite{loh2012high}, the independent case was considered, and the covariance matrix estimator was not analyzed directly, but as part of  a linear regression study. 
In the second case, when the missing data distribution is MCAR but  unknown, we compute estimates of the observation probabilities which can be used to construct an unbiased estimator. This estimator is  the more natural generalization of the sample covariance matrix, and to the best of our knowledge, there are no available results for its estimation error.  A similar approach, however rendering biased  consistent estimators,  has been used in  \cite{kolar2012consistent, cai2016minimax} as  a component of other procedures. 

We also consider a new missing data mechanism, which we call  {\em conditionally MCAR (CMCAR)}.
For  this case, coordinates are observed according to dependent Bernoulli  random variables, but their distributions can change over time,  and  depend  on previous observations.  
A simplified version of this mechanism was proposed in \cite{pavez2018active}, however they used a different estimator, and no error bounds were provided.
The CMCAR mechanism is closely related to active learning, which are data adaptive procedures   designed to maximize system performance, while complying  with problem specific cost constraints.  This mechanism can model situations occurring in  sensor networks, where data dependent missing observations may occur if the signal being measured can affect the sensing process. For example, temperature and humidity may damage a sensor over time.    In other applications,  if an observer is seeking to improve system performance, the subset of variables that are acquired/missed may depend on the data.
In some scenarios,  data acquisition carries a cost and obtaining complete observations might not be possible. For example, in clinical studies, measuring certain features (e.g. acquiring samples or running tests)  may be expensive  \cite{veeramachaneni2005active,lim2005hybrid}, or in sensor networks,   obtaining measurements incurs in  energy consumption  \cite{asif2016matrix, deshpande2004model,gershenfeld2010intelligent}. 
Active feature acquisition techniques have been applied to reduce the cost of   predictive and classification  algorithms while maintaining performance \cite{zheng2002active,saar2009active,melville2005economical,chakraborty2013active}.  
\subsection{Covariance estimation with missing data}

The MCAR mechanism has been considered in previous statistical studies involving the covariance matrix \cite{loh2012high,lounici2014high}.
Particularly in \cite{lounici2014high}, the author considers MCAR observations  with uniform Bernoulli distribution, i.e., all variables are observed with probability $p>0$, and the estimator of (\ref{eq_est_cov_y_mu}). That study shows that the error bound and sample complexity depend on $\erank(\mathbf{\Sigma})/p^2$, thus effectively requiring $1/p^2$ more samples than the sample covariance matrix to achieve the same performance. That work also proposes a regularized estimator for low rank covariance matrices, and shows that by carefully choosing a regularization parameter, its estimation error attains optimal error rates.
A more recent work \cite{park2019non} obtained concentration bounds for the error measured in operator norm when the observations are characterized by dependent Bernoulli random variables with arbitrary observation probabilities. Altough that work is the first to consider the non uniform dependent setting, their error bounds have sub-optimal dependencies on the ambient dimension (an additional linear factor), and the fraction of missing entries (see Table \ref{tab_comparison} and Section \ref{sec:miss_estimator} for more details).

 Cai and Zhang \cite{cai2016minimax} studied estimation of sparse and bandable covariance matrices  under a MCAR missing data mechanism. 
 The obtained error bounds match those for complete observations, up to scale factors that reflect the effect of   missing data. Additionally, the error bounds are pessimistic,  since they depend only on the least observed entry of the covariance matrix. 
Neither  study \cite{lounici2014high,cai2016minimax} offers insight into the relations among the missing data mechanism, the entries of the population covariance matrix and the covariance estimation error. This is because   the missing data mechanism is too simple \cite{lounici2014high}, or the error bounds are too conservative \cite{cai2016minimax}.  
Other works have focused on estimation of sparse inverse covariance matrices \cite{loh2012high,kolar2012consistent}. In particular \cite{loh2012high} derives an error bound for the covariance matrix that is valid for sparse Gaussian graphical models. Missing data in principal component analysis was considered in \cite{lounici2013sparse}.
Table \ref{tab_comparison} lists the most relevant covariance estimators and compares them in terms of sample complexity, that is,  the minimum number of observations required to achieve error below $\epsilon$. 
	\begin{table*}[h]
			\centering
			\scalebox{0.95}{
		\begin{tabular}{|l|l|l|l|l | l|}
			\hline
			Paper	& Structure & Estimator & Missing data  & Sample complexity (uniform case)	& Bound  \\ \hline
			Bunea \& Xiao \cite{bunea2015sample}	& any		& sample covariance	& N/A					&		$\erank(\mathbf{\Sigma}) \log(n)/\epsilon^2$		&  E 	\\ \hline 
			Cai \& Zhang \cite{cai2016minimax}	& bandable	& block threshold	& MCAR		&	$(1/\widehat{p}_{\min}) \left( C\log(n)/\epsilon^2 \vee (C/\epsilon^2)^{1+1/2\alpha} \right)$				& E\\
			Cai \& Zhang \cite{cai2016minimax}	& sparse		& adaptive threshold	& MCAR	&	$ \rho_{n,N}\log(n)/\epsilon^2 $				&E\\
			Lounici \cite{lounici2014high}& any		& Eq.(\ref{eq_est_cov_y_mu})	& UI MCAR	&	$ ( {\erank(\mathbf{\Sigma})}/{p^2})\log(2n)\left( {1}/{\epsilon}  \vee (\log(n)+\log(N)+c_1 p)\right)^2$				& P \\ 
			Park \& Lim \cite{park2019non}& any & Eq.(\ref{eq_est_cov_y_mu})	& NUD MCAR& $( {\erank(\mathbf{\Sigma})}/{p^{1/2}})\log(2n)\left( {n}/{p^{3/2}\epsilon}  \vee (\log(n)+\log(N)+c_1 p)\right)^2$ & P \\
			\hline 
			Theorem \ref{th_expected_error}&any		&Eq.(\ref{eq_est_cov_y_mu})					&NUD MCAR	& 	$ (\erank( \mathbf{\Sigma})/ p^2) \log(n) (e\log(N) \vee 2C_2/\epsilon)^2$ &E \\ \hline 
			Theorem \ref{th_expected_error_MAR} & any & Eq.(\ref{eq_est_cov_y_mu_mcar2}) &  MCAR & ${\srank}_{\min}(\mathbf{\Sigma, \widehat{P}}) \log(n) (e\rho \log(N) \vee 2C_2/\epsilon)^2 $ & E \\ \hline 
		\end{tabular}
}
		\caption{Comparison of   covariance estimators with missing data.  For Cai \& Zhang \cite{cai2016minimax}, the constants  $\alpha$ and $\rho_{n,N}$  parameterize the classes of bandable and sparse covariance matrices respectively. $\widehat{p}_{\min}$ is the smallest non zero entry of $\widehat{\mathbf{ P}}$. Some universal numerical constants have been ignored, for detailed results see corresponding Theorems in Section \ref{sec:miss_missing}. Abbreviations:    uniform independent (UI),  non uniform dependent (NUD), expectation bound (E), high probability bound (P).  }
		\label{tab_comparison}
		 \end{table*}
A particular version of the CMCAR  mechanism and a covariance estimator was introduced in  \cite{pavez2018active},  but no error bounds were provided.
\section{Covariance estimators}
 \label{sec:miss_estimator}
%
\subsection{Notation}
We denote matrices and vectors by bold letters. For a matrix $\mathbf{X} = (x_{ ij} )$ its entries are denoted
by $x_{ij}$ or $\mathbf{X}_{ij}$ and for a vector $\mathbf{x}$ its i-th entry is denoted by $x_i$. The entry-wise product\footnote{The entry-wise product is also known as Hadamard or Schur product.} between matrices is defined as $(\mathbf{A} \odot \mathbf{B})_{ij} = a_{ij} b_{ij}$. We use $\Vert \cdot \Vert_q$  for entry-wise matrix norms, with  $q=2$ corresponding to  the Frobenius norm. $\Vert \cdot \Vert$ denotes the $\ell_2$ norm when applied to vectors, and the operator  norm (largest singular value) when applied to matrices.  The nuclear norm (sum of singular values) is denoted by $\Vert \cdot \Vert_{\star}$.
The set of integers $\lbrace 1, 2, \cdots, d \rbrace$ is abbreviated by $[d]$.
\subsection{Missing observations}
The complete data is represented with the  $n \times N$  matrix $\mathbf{X}$. Its $i$-th row and $k$-th column are the vectors  $\mathbf{x}_i$ and $\mathbf{x}^{(k)}$ respectively.  
The $n \times N$ matrices  $\mathbf{Y}$ and  $\mathbf{\Delta}$ denote the observation and missing data patterns respectively. $\mathbf{\Delta}$ has entries in $\lbrace 0, 1 \rbrace$.   The $k$th columns of $\mathbf{Y}$ and  $\mathbf{\Delta}$  are denoted by the vectors  $\mathbf{y}^{(k)}$ and $\boldsymbol\delta^{(k)}$ respectively. These quantities are related  through
\begin{equation*}
\mathbf{y}^{(k)} = \boldsymbol\delta^{(k)} \odot \mathbf{x}^{(k)}.
\end{equation*}
When the $j$-th entry of $\boldsymbol\delta^{(k)}$ is equal to $0$, the corresponding entry of $\mathbf{x}^{(k)}$ is missing. 
If the column vectors $\mathbf{x}^{(k)}$ are i.i.d., zero mean and with population covariance matrix $\mathbf{\Sigma}$, the sample covariance matrix
of the complete data 
\begin{equation*}
\mathbf{S} = \frac{1}{N} \mathbf{ X X^{\top}} = \frac{1}{N}\sum_{k=1}^N \mathbf{x}^{(k)} {\mathbf{x}^{(k)}}^{\top},
\end{equation*}
is an unbiased estimator for $\mathbf{\Sigma}$. The  behavior of $\mathbf{S}$ for finite $N$  has been studied extensively and it is well understood \cite{lounici2014high,bunea2015sample,koltchinskii2017concentration}. In the rest of this section we will introduce alternative  estimators for $\mathbf{\Sigma}$ that depend on  the observations  $\mathbf{Y}$, and properties of the  missing data pattern $\mathbf{\Delta}$.
\subsection{Missing data mechanisms and covariance estimators}
Consider a random vector  $\mathbf{x}=[x_1,x_2,\cdots,x_n]^{\top}$  in $\mathbb{R}^n$ with population covariance matrix $\mathbf{\Sigma}$. The  purpose of a covariance estimator is to  accurately  approximate $\mathbf{\Sigma}$  from $\lbrace \mathbf{x}^{(k)} \rbrace^{N}_{k=1}$, a set of independent identically distributed (i.i.d.) copies of $\mathbf{x}$. 
Here, however, we only have access to the sequence of incomplete observations $\lbrace \mathbf{y}^{(k)} \rbrace^{N}_{k=1}$.
Our interest is in estimators with the following desirable properties:
\begin{enumerate}
	\item {\em Unbiased:} their expected value is equal to  $\mathbf{\Sigma}$.
	\item {\em Consistent:} as the number of observations $N\rightarrow \infty$ the estimator converges to $\mathbf{\Sigma}$.
	\item {\em Minimum sample complexity:}  minimum number of samples to achieve error below $\epsilon>0$.
\end{enumerate}
We  consider two missing data mechanisms, and  three different estimators which are described next. 
%
%
%
%
%
\paragraph{\NEW{Unbiased  estimator under MCAR observations with known observation probabilities}}
The MCAR model asserts that $\mathbf{X}$ and $\mathbf{\Delta}$ are independent \cite{little2014statistical}.
We model the missing data  with a random vector $ \boldsymbol\delta = (\delta_i)$ of dependent Bernoulli $0-1$ random variables, where $\Prob(\delta_i =1, \delta_j =1)=\E[\delta_i \delta_j ]=p_{ij}$. 
%
%
We assume that $\boldsymbol\delta^{(k)}$ are i.i.d. copies of $\boldsymbol\delta$. 
If all $p_{ii} = 1$, then $\mathbf{y}^{(k)} = \mathbf{x}^{(k)}$ for all $k$, thus we have perfect observation of $\mathbf{x}^{(k)}$.  We are interested in recovering the covariance matrix of  $\mathbf{x}$ from $\lbrace \mathbf{y}^{(k)} \rbrace_{k=1}^N$, when $0<p_{ij} \leq 1$   for all $i  \in [n]$ and $j  \in [n]$.  
Let $\mathbf{P} = (p_{ij})$, and $\mathbf{\Gamma}$ its Hadamard (entry-wise) inverse, that is $\mathbf{P} \odot \mathbf{\Gamma} = \mathbf{1} \mathbf{1}^{\top}$, which is guaranteed to exist since $p_{ij}>0$ for all $i,j \in [n]$.  When $\E[\mathbf{x}] = \mathbf{0}$, we  will use
\begin{equation}\label{eq_est_cov_y_mu}
\widehat{\mathbf{\Sigma}} = \frac{1}{N} \sum_{k=1}^N \mathbf{y}^{(k)}  {\mathbf{y}^{(k)}}^{\top} \odot \mathbf{\Gamma}.
\end{equation}
Next we show that (\ref{eq_est_cov_y_mu}) is unbiased. We consider a slight generalization and assume the population mean is known and non-zero given by $\E[\mathbf{x}] = \boldsymbol\mu$. The estimator (\ref{eq_est_cov_y_mu}) becomes
		\begin{equation*}
		\widehat{\mathbf{\Sigma}}=\frac{1}{N} \sum_{k=1}^N (\mathbf{y}^{(k)} - \boldsymbol\delta^{(k)} \odot \boldsymbol\mu )(  {\mathbf{y}^{(k)}}- \boldsymbol\delta^{(k)} \odot \boldsymbol\mu )^{\top} \odot \mathbf{\Gamma}.
		\end{equation*}
The centered variable satisfies $\mathbf{y}^{(k)} - \boldsymbol\delta^{(k)} \odot \boldsymbol\mu =  \boldsymbol\delta^{(k)} \odot(\mathbf{x}^{(k)} -  \boldsymbol\mu) $, so we may write its expectation
	\begin{align*}
	&\E[(\mathbf{y}^{(k)} - \boldsymbol\delta^{(k)} \odot \boldsymbol\mu )(  {\mathbf{y}^{(k)}}- \boldsymbol\delta^{(k)} \odot \boldsymbol\mu )^{\top} ] \\
	&=\E[(\boldsymbol\delta^{(k)} \odot(\mathbf{x}^{(k)} -  \boldsymbol\mu)  )(\boldsymbol\delta^{(k)} \odot(\mathbf{x}^{(k)} -  \boldsymbol\mu)  )^{\top} ] \\
	&=\E[\boldsymbol\delta^{(k)}{\boldsymbol\delta^{(k)}}^{\top} \odot (\mathbf{x}^{(k)} -  \boldsymbol\mu)  (\mathbf{x}^{(k)} -  \boldsymbol\mu)  ^{\top} ] \\
	&=\E[\boldsymbol\delta^{(k)}{\boldsymbol\delta^{(k)}}^{\top}] \odot \E[(\mathbf{x}^{(k)} -  \boldsymbol\mu)  (\mathbf{x}^{(k)} -  \boldsymbol\mu)  ^{\top}] = \mathbf{P} \odot \mathbf{\Sigma}.
	\end{align*}
	With this we can show that  expectation of $\widehat{\mathbf{\Sigma}}$ is $\mathbf{\Sigma}$, hence we have the following.
\begin{proposition}\label{prop_unbiased_MCAR1}
	 (\ref{eq_est_cov_y_mu}) is an unbiased estimator for $\mathbf{\Sigma}$ .
\end{proposition}
Our error bound for (\ref{eq_est_cov_y_mu}) is given in Theorem \ref{th_expected_error}.
Previously, 
\cite{lounici2014high} characterized the estimation error of (\ref{eq_est_cov_y_mu})   when $\boldsymbol\delta$ has independent and identically distributed coordinates. Later,  \cite{park2019non} derived error bounds for the non uniform dependent setting, however, as we will see in Section \ref{sec:miss_missing}, their bounds do not match those of \cite{lounici2014high} or ours, and thus are sub-optimal.
\paragraph{\NEW{Unbiased  estimator under MCAR observations with unknown observation probabilities}}
Here we consider a more realistic scenario in which the observation probabilities are unknown and have to be estimated from data. For this we define the matrix of empirical observation probability 
\begin{equation}
\widehat{\mathbf{P}} = \frac{1}{N}  \sum_{k=1}^N \boldsymbol{\delta}^{(k)}{\boldsymbol{\delta}^{(k)}}^{\top}, 
\end{equation}
and   denote its entries by $\hat{p}_{ij}$. The quantity $\hat{p}_{ij}$   is equal to zero when the $i$-th and $j$-th variables are not observed together. Since  we can only  estimate  the entries of the covariance matrix for which there are observations,  we define the set  $\M = \lbrace (i,j): \hat{p}_{ij} \neq 0 \rbrace$. \NEW{Let $\mathbf{\Sigma}_{\M}$ be matrix whose entries in $\M$ coincide with the population covariance,  and are zero otherwise.} Let  $\widehat{\mathbf{\Gamma}} = (\hat{\gamma}_{ij}) $ be the Hadamard inverse of $\max(\widehat{\mathbf{P}},\frac{1}{N}\mathbf{11^{\top}})$, which satisfies $\hat{\gamma}_{ij} \hat{p}_{ij} = 1$ if $(i,j) \in \M$, and  $\hat{\gamma}_{ij} \hat{p}_{ij} = 0$ otherwise. An  unbiased estimator for $\mathbf{\Sigma}_{\M}$ is  given by
\begin{equation}\label{eq_est_cov_y_mu_mcar2}
\widehat{\mathbf{\Sigma}} = \frac{1}{N} \sum_{k=1}^N \mathbf{y}^{(k)}  {\mathbf{y}^{(k)}}^{\top} \odot \widehat{\mathbf{\Gamma}}.
\end{equation}
In fact, we can use the more general estimator
\begin{equation*}
		\widehat{\mathbf{\Sigma}}=\frac{1}{N} \sum_{k=1}^N (\mathbf{y}^{(k)} - \boldsymbol\delta^{(k)} \odot \boldsymbol\mu )(  {\mathbf{y}^{(k)}}- \boldsymbol\delta^{(k)} \odot \boldsymbol\mu )^{\top} \odot \widehat{\mathbf{\Gamma}}
\end{equation*}
 when the  population mean $\boldsymbol\mu$ is non zero. 
To see that this estimator is unbiased, we write each term in the sum as
$(\mathbf{x}^{(k)} - \boldsymbol\mu )(  {\mathbf{x}^{(k)}}-  \boldsymbol\mu )^{\top} \odot  \boldsymbol\delta^{(k)}{ \boldsymbol\delta^{(k)}}^{\top} \odot \widehat{\mathbf{\Gamma}}$. Applying conditional expectation,  and independence of $\mathbf\Delta$ and $\mathbf{X}$ (MCAR assumption) we get
\begin{align*}
			&\E[(\mathbf{x}^{(k)} - \boldsymbol\mu )(  {\mathbf{x}^{(k)}}-  \boldsymbol\mu )^{\top} \odot  \boldsymbol\delta^{(k)}{ \boldsymbol\delta^{(k)}}^{\top} \odot \widehat{\mathbf{\Gamma}} \mid \mathbf{\Delta }] \\
			&=\E[(\mathbf{x}^{(k)} - \boldsymbol\mu )(  {\mathbf{x}^{(k)}}-  \boldsymbol\mu )^{\top} \mid \mathbf{\Delta }]\odot  \boldsymbol\delta^{(k)}{ \boldsymbol\delta^{(k)}}^{\top} \odot \widehat{\mathbf{\Gamma}} \\
			&= \mathbf{\Sigma}\odot  \boldsymbol\delta^{(k)}{ \boldsymbol\delta^{(k)}}^{\top} \odot \widehat{\mathbf{\Gamma}}.
\end{align*}
Adding up all terms produces
\begin{equation*}
			\E[	\widehat{\mathbf{\Sigma}}\mid\mathbf{\Delta}] = \frac{1}{N} \sum_{k=1}^N  \mathbf{\Sigma} \odot \boldsymbol\delta^{(k)}{ \boldsymbol\delta^{(k)}}^{\top} \odot \widehat{\mathbf{\Gamma}} =  \mathbf{\Sigma} \odot  \widehat{\mathbf{P}}\odot \widehat{\mathbf{\Gamma}} = \mathbf{\Sigma}_{\M}.
\end{equation*}
If   all  $\hat{p}_{ij}>0$, then $\widehat{\mathbf{P}}\odot \widehat{\mathbf{\Gamma}} = \mathbf{1} \mathbf{1}^{\top}$,  $\mathbf{\Sigma}_{\M} = \mathbf{\Sigma}$ and $\E[\widehat{\mathbf{\Sigma}}] = \mathbf{\Sigma}$.
%
\begin{proposition}\label{prop_unbiased_MCAR2}
	Conditioned on $\mathbf{\Delta}$,  (\ref{eq_est_cov_y_mu_mcar2}) is an unbiased estimator for $\mathbf{\Sigma}_{\M}$. If all  $\hat{p}_{ij}>0$, then (\ref{eq_est_cov_y_mu_mcar2}) is unbiased for $\mathbf{\Sigma}$.
\end{proposition}
A similar  but biased estimator that does not know the missing data distribution was used as  part of the banded and sparse covariance estimators from \cite{cai2016minimax}, and the inverse covariance estimator from \cite{kolar2012consistent}.  However, to the best of our knowledge we are the first to establish estimation error bounds  for an estimator of this form in Theorem \ref{th_expected_error_MAR}.
\paragraph{Unbiased estimator under CMCAR observations}
The CMCAR mechanism presented here generalizes the  MCAR to allow  missing data with a time varying distribution. One instance of this mechanism are active sequential observations, where the observation probabilities depend on previously observed data. Another,  is a MCAR  mechanism with time varying observation probabilities.  To the best of our knowledge, this type of missing data mechanism has not been considered for  covariance matrix estimation.

The index $t=1,2, \cdots$ denotes time  (or iteration number).   The observed vectors at time $t$ are $\lbrace \mathbf{y}^{(t,1)}, \cdots, \mathbf{y}^{(t,N_t)} \rbrace$, which can be stored as columns of the  $n \times N_t$ dimesional matrix $\mathbf{Y}_t$. Each of these vectors is defined as $\mathbf{y}^{(t,k)} = \boldsymbol{\delta}^{(t,k)} \odot \mathbf{x}^{(t,k)}$, where  $\left\lbrace\boldsymbol{\delta}^{(t,k)} \right\rbrace_{k=1}^{N_t}$ are  identically distributed Bernoulli with sampling probabilities ${p}_{ij}^{(t)}$.   
 We assume that for $t \geq 1$, 
\begin{enumerate}
	\item {\em Complete data is i.i.d.}. For $k \in [N_t]$, and $t \in [T]$, the complete observations $\mathbf{x}^{(t,k)}$ are i.i.d. copies of $\mathbf{x}$.
	\item {\em Stochastic observation parameters. } The sampling probabilities are a function of previous observations, thus ${p}_{ij}^{(t)} = f_{ijt}(\mathbf{Y}_{t-1}, \mathbf{Y}_{t-2},\cdots, \mathbf{Y}_1)$.
	\item {\em Conditionally i.i.d. observations.} Conditioned on the event $\lbrace \mathbf{Y}_{\tau}: \tau <t \rbrace$, the random vectors $\mathbf{y}^{(t,1)}, \cdots, \mathbf{y}^{(t,N_t)}$ are i.i.d..
	\item {\em Conditionally MCAR  (CMCAR) mechanism.} Conditioned on the event $\lbrace \mathbf{Y}_{\tau}: \tau <t \rbrace$, the  vectors $\left\lbrace\boldsymbol{\delta}^{(t,k)} \right\rbrace_{k=1}^{N_t}$ are i.i.d.,  and  independent of $\mathbf{x}^{(t,k)}$.
\end{enumerate} 
%
%
The sampling probabilities are positive (i.e., $p_{ij}^{(t)}>0$ for all $i,j,t$). We will also use the  matrices $\mathbf{P}_{t}=(p_{ij}^{(t)})$, and $\mathbf{\Gamma}_{t} $, which are Hadamard inverses, that is,  $\mathbf{P}_{t} \odot \mathbf{\Gamma}_{t} = \mathbf{1} \mathbf{1}^{\top}$.
For simplicity, conditional expectation is abbreviated as
\begin{equation*}
\mathbb{E}_{t-1}[~\cdot~] = \mathbb{E}[~\cdot~ \mid \mathbf{Y}_{\tau}: \tau <t].
\end{equation*}
Below we present some useful quantities
\begin{align}\label{eq_cov_est_local}
\mathbf{Z}_t 				&=  \mathbf{Y}_t \mathbf{Y}_t^{\top} \odot \mathbf{\Gamma}_t  \\
\mathbf{W}_T 				&= \sum_{t=1}^T (\mathbf{Z}_t - N_t\mathbf{\Sigma}) \nonumber\\
\widehat{\mathbf{\Sigma}}_T	&= \frac{\sum_{t=1}^T \mathbf{Z}_t}{\sum_{t=1}^T N_t}. \label{eq_cov_est_martinga}
\end{align}
Some facts about $\mathbf{Z}_t,~ \mathbf{W}_t$ and  $\widehat{\mathbf{\Sigma}}_T $ are stated below.
\begin{enumerate}
	\item {\em $\mathbf{Z}_t$ is an unbiased estimator for $N_t\mathbf{\Sigma}$}. This follows from an application of  conditional expectation, and using  (\ref{eq_est_cov_y_mu}), 
	\begin{align}
	{\E}_{t-1}[\mathbf{Z}_t] &= {\E}_{t-1}\left[ \mathbf{Y}_t \mathbf{Y}_t^{\top} \odot \mathbf{\Gamma}_t \right] \nonumber\\
	&=  {\E}_{t-1}\left[ \mathbf{Y}_t \mathbf{Y}_t^{\top}\right] \odot \mathbf{\Gamma}_t \nonumber\\
	&= N_t \mathbf{\Sigma} \odot \mathbf{P}_t \odot \mathbf{\Gamma}_t \nonumber\\
	&= N_t \mathbf{\Sigma}. \label{eq_martinga_cond_expect_zero}
	\end{align}
	\item {\em $\mathbf{W}_T$ is a zero mean matrix martingale}. The martingale property states that $\mathbb{E}_T[\mathbf{W}_{T+1}] = \mathbf{W}_T$. We can verify this applying conditional expectation,
	\begin{equation*}
	\mathbb{E}_T[\mathbf{W}_{T+1}] = {\E}_{T}[\mathbf{Z}_{T+1} - N_{T+1}\mathbf{\Sigma}] + {\E}_{T}[\mathbf{W}_T] = \mathbf{W}_T.
	\end{equation*}
This property implies  that $\mathbb{E}[\mathbf{W}_{T}] =\mathbb{E}[\mathbf{W}_{1}] =\mathbf{0}$.
	\item {\em $\widehat{\mathbf{\Sigma}}_T$ is an unbiased estimator for $\mathbf{\Sigma}$}. This is a direct consequence of   property 2) and the identity ${\mathbf{W}}_T = (\widehat{\mathbf{\Sigma}}_T - \mathbf{ \Sigma}) (\sum_{t=1}^T N_t)$.
\end{enumerate}
Performance of this estimator is presented in Theorem \ref{th_expected_error_martin}.
\NEW{
Given Propositions \ref{prop_unbiased_MCAR1} and \ref{prop_unbiased_MCAR2}, the observations $\mathbf{y}^{(k)}$ can always be centered to have mean zero, provided the population mean is known. For the rest of the paper we will assume, without loss of generality,  that the mean is known and $\boldsymbol\mu = \mathbf{0}$.  
}
\NEW{
\subsection{Unknown  mean}
When the mean is unknown, centering using an estimate for the population mean is not sufficient to construct an unbiased estimator. To see this, consider 
\begin{equation}\label{eq_estimator_cai_kolar}
\frac{1}{N}\sum_{k=1}^N (\mathbf{y}^{(k)} - \boldsymbol\delta^{(k)}\odot \bar{\mathbf{y}} )(\mathbf{y}^{(k)} - \boldsymbol\delta^{(k)}\odot \bar{\mathbf{y}} )^{\top } \odot \widehat{\mathbf{\Gamma}},
\end{equation}
which was utilized in  \cite{cai2016minimax,kolar2012consistent}  when the mean and missing data distributions are unknown.
Although (\ref{eq_estimator_cai_kolar}) is a consistent estimator for $\mathbf{\Sigma}$  \cite[Lemma 2.1]{cai2016minimax}, it is also biased. To see the latter, notice that in the complete data case, (\ref{eq_estimator_cai_kolar}) becomes $({1}/{N})\sum_{k=1}^N (\mathbf{x}^{(k)} -  \bar{\mathbf{x}} )(\mathbf{x}^{(k)} -  \bar{\mathbf{x}} )^{\top }$, which has expectation equal to $({N-1}/{N}) \mathbf{\Sigma}$. A similar situation is encountered when $\mathbf{P}$ is known. In  appendix \ref{sec:miss_non_zero_mean}, we  discuss unbiased alternatives to (\ref{eq_estimator_cai_kolar}),  for known and unknown  missing data distribution. However, the expressions are much more cumbersome,  thus  a theoretical error study is left for future work. 
}
%
%
%
%
%
%
%
%
%
%
\section{Error bounds}
 \label{sec:miss_missing}
We start this section by  introducing some concepts related to sub-Gaussian vectors and notions of complexity of covariance matrices. In Theorems \ref{th_expected_error}, \ref{th_expected_error_MAR}, and \ref{th_expected_error_martin} we present error bounds  in operator norm for each estimator.    All   results of this section  depend on new quantities we denote {\em scaled effective rank}, which generalize the effective rank to account for partial observations. All our bounds can be viewed as  non trivial applications of a moment inequality found in \cite{chen2012masked}.  After  some simplifications, our results are summarized in  Table \ref{tab_comparison}.
Based on  Propositions \ref{prop_unbiased_MCAR1} and \ref{prop_unbiased_MCAR2},  in this section we assume the population mean is known and equal to zero. 
\subsection{Preliminaries and  assumptions} 
We first introduce a proxy for the rank of a matrix.
\begin{definition}\label{def_erank} The {\bf effective rank} \cite{lounici2014high,vershynin2016high} of a $n \times n$ matrix $\mathbf{A}$ is defined as
	\begin{equation*}
	{\erank}(\mathbf{A}) = \frac{\Vert \mathbf{A} \Vert_{\star}}{\Vert \mathbf{A}\Vert}.
	\end{equation*}
\end{definition}
The effective rank obeys $1 \leq \erank(\mathbf{A}) \leq \rank(\mathbf{A}) \leq n$.  
When $\mathbf{A}$ is positive semi-definite, $\Vert \mathbf{A}\Vert_{\star}=\tr(\mathbf{A})$.  The effective rank, as we will see in latter sections,  is the  parameter that quantifies the sample complexity of covariance estimation.

We introduce sub-Gaussian variables and vectors  below. 
\begin{definition}[\cite{vershynin2016high}]
	If $\E[\exp(z^{2} /K^{2})] \leq 2$ holds for  some $K>0$, we say $z$ is sub-Gaussian. If  $\E[\exp(\vert z \vert /K)] \leq 2$ holds for some $K>0$, we say $z$ is sub-exponential.
\end{definition} 
\begin{definition}[\cite{vershynin2016high}] For a real valued random variable $z$, the sub-Gaussian and sub-exponential norms are defined as
	\begin{equation*}
	\Vert z \Vert_{\Psi_{\alpha}} = \inf\lbrace u>0 : \E[\exp(\vert z\vert^{\alpha}/u^{\alpha})] \leq 2
	\rbrace,
	\end{equation*}
	for $\alpha =2$ and $\alpha =1$ respectively.
\end{definition}

\begin{definition}
	If a  random vector $\mathbf{z}$ in $\mathbb{R}^n$ satisfies
	\begin{equation*}
	\Vert \mathbf{z} \Vert_{\Psi_{\alpha}} = \sup_{\mathbf{u}: \Vert \mathbf{u}\Vert=1} \Vert  \mathbf{u}^{\top} \mathbf{z}\Vert_{\Psi_{\alpha}} <\infty,
	\end{equation*}
	it is called sub-Gaussian or sub-exponential, for $\alpha=2$ or $\alpha =1$ respectively.
\end{definition}
We will work under Assumption \ref{assum_trace} detailed below, which is standard in the covariance estimation literature  \cite{bunea2015sample,lounici2014high}.
\begin{assumption}\label{assum_trace}
	The distribution of $\mathbf{ x}$ obeys.
	\begin{itemize}
		\item $\mathbf{x}$ is a sub-Gaussian vector, therefore   there is a constant $c$ so that $\mathbb{E}[\vert \mathbf{u}^{\top} \mathbf{x}\vert^r]^{1/r} \leq c \Vert  \mathbf{u}^{\top} \mathbf{x}\Vert_{\Psi_{2}} \sqrt{ r} $ for every  unit vector $\mathbf{ u} \in \mathbb{R}^n$, and $r\geq 1$.
		\item There is a constant $c_0$ so that for all $\mathbf{u} \in \mathbb{R}^n$, then $ \Vert  \mathbf{u}^{\top} \mathbf{x}\Vert_{\Psi_{2}}^2 c_0 \leq  \mathbb{E}(\mathbf{u}^{\top} \mathbf{x})^2$
	\end{itemize}
\end{assumption}
%
%
%
%
Finally we introduce  a generalization of  the   effective rank.
\begin{definition} A {\bf scaled effective rank} of a matrix  $\mathbf{\Sigma}$, is any quantity that takes values in $[\erank(\mathbf{\Sigma}), \erank(\mathbf{\Sigma})/\rho]$ for some constant $0<\rho<1$, where the constant does not depend on $\mathbf{\Sigma}$.
\end{definition}
Examples that will appear in this section are 
\begin{align*}
&{\srank}_{\min}(\mathbf{\Sigma, P}) = \frac{1}{\Vert \mathbf{\Sigma} \Vert} \max_{j \in [n]} \sum_{i=1}^n \frac{\Sigma_{ii}}{p_{ij}}=\frac{1}{\Vert \mathbf{\Sigma} \Vert}\Vert \diag(\mathbf{\Sigma})\mathbf{\Gamma} \Vert_{1 \rightarrow 1},\\
&{\srank}_2(\mathbf{\Sigma, P})  = \frac{1}{\Vert \mathbf{\Sigma} \Vert}\Vert \diag(\mathbf{\Sigma})^{1/2} \mathbf{\Gamma}\diag(\mathbf{\Sigma})^{1/2}  \Vert_F,\\
&\tilde{\srank}_{\min}(\mathbf{\Sigma, p}) = \frac{1}{\Vert \mathbf{\Sigma} \Vert} \frac{1}{\min_i p_{ii}} \sum_{i=1}^n \frac{\Sigma_{ii}}{p_{ii}},\\
&\tilde{\srank}_2(\mathbf{\Sigma, p})  = \frac{1}{\Vert \mathbf{\Sigma} \Vert} \sum_{i=1}^n \frac{\Sigma_{ii}}{p_{ii}^2},
\end{align*}
where $\mathbf{p}$ is a vector containing the diagonal entries of $\mathbf{P}$.
%
\subsection{\NEW{MCAR case with known distribution}}
Now we state our main result for the MCAR model.
\begin{theorem}\label{th_expected_error}
	Suppose the random vector $\mathbf{x}$ obeys Assumption \ref{assum_trace}. If all $p_{ij}>0$, 
	and  also  $N \geq e^2$ and $n \geq 3$, we have that (\ref{eq_est_cov_y_mu}) satisfies
	\begin{align}\label{spectral_bound}
	&\mathbb{E}\left[ \Vert\widehat{\mathbf{\Sigma}} - \mathbf{\Sigma}\Vert^2 \right]^{1/2} \leq C \Vert\mathbf{\Sigma}\Vert  \\
	&\left( \sqrt{\frac{8e\log(n) {{\srank}_{\min}}(\mathbf{\Sigma,P}) }{N}} + \frac{4e^2\log(n)\log(N){\srank_2}(\mathbf{\Sigma,P})}{N}\right), \nonumber
	\end{align}
	where $C$ is an  universal constant.
	In addition, if $\boldsymbol\delta$ has independent entries, the parameters ${{\srank}_{\min}}(\mathbf{\Sigma,P})$ and ${\srank_2}(\mathbf{\Sigma,P})$ in  (\ref{spectral_bound}) can be replaced by $\tilde{\srank}_{\min}(\mathbf{\Sigma,p})$ and $\tilde{\srank}_2(\mathbf{\Sigma,p})$, respectively.
\end{theorem}
\NEW{
The upper bound from Theorem \ref{th_expected_error} suggests that if variable $x_i$ has high   variance (high $\Sigma_{ii}$) or if it is strongly correlated with other variables (high $\Sigma_{ij}$), then the effect of missing data can be more severe, since the upper bound increases with $\Sigma_{ii}/p_i$ and $\Sigma_{ij}/p_{ij}$. 
}

 Theorem \ref{th_expected_error} is an improvement  over previous bounds under  the MCAR mechanism with known missing data distribution. In particular, we obtain better convergence rates than \cite{park2019non}, and consider a more general missing data mechanism than \cite{lounici2014high}. However, these previous studies obtained high probability  error bounds, while here we bound the expectation. 
\paragraph{ Comparison with  covariance estimation with complete observations} 
In this case,  $p_{ii}=1$ for all $i \in [n]$, which implies  that all $p_{ij}=1$, and the scaled effective rank coincides with the effective rank.  Several recent papers \cite{lounici2014high,bunea2015sample,koltchinskii2017concentration} have obtained bounds in operator norm that scale as
\begin{equation*}
\mathcal{O}\left(\sqrt{\frac{\erank(\mathbf{\Sigma})}{N}} + \frac{\erank(\mathbf{\Sigma})}{N}\right)=\mathcal{O}\left(\sqrt{\frac{\erank(\mathbf{\Sigma})}{N}} \vee \frac{\erank(\mathbf{\Sigma})}{N}\right),
\end{equation*}
where the $\mathcal{O}$ notation hides constants and logarithmic dependencies on $n$ and $N$. 
Our bound in spectral norm matches those  of \cite{lounici2014high,bunea2015sample} up to  $\log(N)$ factors. For the Gaussian case, it was shown in  \cite{koltchinskii2017concentration} that dimension free bounds can be attained (without any $\log$ terms). It is also worth mentioning that the  error bounds from \cite{lounici2014high,bunea2015sample} are high probability results based on variations of matrix Bernstein inequalities for the operator norm. In \cite{lounici2014high}, one is used for sums of bounded random matrices, while   \cite{bunea2015sample} applies a version  for unbounded matrices. Our results  are based in a moment inequality for sums of independent, possibly unbounded matrices from \cite{chen2012masked}, thus significantly simplifying the proofs.
\paragraph{ Sample complexity}
We can guarantee covariance estimation error in operator norm of magnitude
\begin{equation*}
\mathbb{E}\left[ \Vert\widehat{\mathbf{\Sigma}} - \mathbf{\Sigma}\Vert^2 \right]^{1/2}  \leq  \epsilon \Vert\mathbf{\Sigma}\Vert,
\end{equation*}
if the number of observations obeys
\begin{equation}\label{eq_mcar_sample_complexity}
{N} \geq 8e\log(n) {{\srank}_{\min}}(\mathbf{\Sigma, P})  \left(e \rho \log(N) \vee \frac{2C_2}{\epsilon}\right)^2,
\end{equation}
where $\rho = {{\srank}_{2}}(\mathbf{\Sigma, P}) / {{\srank}_{\min}}(\mathbf{\Sigma, P})$.
Compared with the case of complete observations, the missing data case requires extra samples by a factor of ${{\srank}_{\min}}(\mathbf{\Sigma, P})/\erank(\mathbf{\Sigma})$  to attain equal accuracy. When $\boldsymbol\delta$ has independent coordinates, a simpler sufficient condition is given by
\begin{equation*}
{N} \geq 8e\log(n) {\tilde{\srank}_{\min}}(\mathbf{\Sigma, p})  \left(e  \log(N) \vee \frac{2C_2}{\epsilon}\right)^2.
\end{equation*}
This condition  can be derived from  the looser bound that depends on $\tilde{\srank}_{2}(\mathbf{\Sigma, p})$ and  $\tilde{\srank}_{\min}(\mathbf{\Sigma, p})$. In this case, both $\tilde{\srank}_{2}(\mathbf{\Sigma, p})$ and ${\srank}_{\min}(\mathbf{\Sigma, p})$ are upper bounded by $\tilde{\srank}_{\min}(\mathbf{\Sigma, p})$.
\paragraph{ Parameter scaling} Each entry of the sample covariance matrix is a (weighted) average of a different number of observations. In particular the $ij$ off-diagonal element of $\widehat{\mathbf{\Sigma}}$ is observed  $Np_{ij}$ times on average. This reduced data acquisition rate is reflected in an increased error bound of the same order.
\paragraph{Comparison with   \cite{lounici2014high} and \cite{park2019non}}
 Since the bounds from \cite{lounici2014high} assume uniform independent observation variables, we will compare under that setting. Lounici \cite{lounici2014high} establishes that with probability at least $1-e^{-\nu}$
 \begin{align}\label{eq_lounici_uniform}
    & \Vert\widehat{\mathbf{\Sigma}} - \mathbf{\Sigma}\Vert  \leq C \Vert \mathbf{\Sigma}\Vert  \\
     & \max \left\lbrace \sqrt{\frac{ {\erank}(\mathbf{\Sigma})(\nu+\log(2n))  }{p^2 N}  }\right. , \nonumber \\
	&\left.  \frac{ {\erank}(\mathbf{\Sigma})(\nu+\log(2n)) (C_1p + \nu +\log(N)) }{p^2 N} \right\rbrace, \nonumber
 \end{align}
 where $C, C_1>0$ are universal constants. On the other hand, Park and Lim  \cite{park2019non} derived a bound under MCAR observations, which after simplification under the uniform and independent assumption is reduced to
 \begin{align}\label{eq_park_uniform}
      & \Vert\widehat{\mathbf{\Sigma}} - \mathbf{\Sigma}\Vert  \leq C \Vert \mathbf{\Sigma}\Vert n \ \\
     & \max \left\lbrace \sqrt{\frac{ {\erank}(\mathbf{\Sigma})(\nu+\log(2n))  }{p^{7/2} N}  }\right. , \nonumber  \\
	&\left.  \frac{ {\erank}(\mathbf{\Sigma})(\nu+\log(2n)) (C_1p + \nu +\log(N)) }{p^2 N} \right\rbrace,\nonumber
 \end{align}
 also with probability at least $1-e^{\nu}$. Setting  $\nu = \log(2n)$, and for large enough $N$, (\ref{eq_lounici_uniform}) and Theorem \ref{th_expected_error} imply that the estimation  error scales as
 \begin{equation*}
     \mathcal{O}\left( \sqrt{\frac{ {\erank}(\mathbf{\Sigma})\log(2n)  }{p^2 N}  } \right),
 \end{equation*}
 while the bound  (\ref{eq_park_uniform}) has a worst rate by a factor of 
 ${n}/{p^{3/4}}$. 
 We can also compare in terms of sample complexity. Lounici's bound implies the following sample complexity condition
 \begin{equation*}
     N \geq C \frac{\erank(\mathbf{\Sigma})}{p^2}\log(2n)\left( \frac{1}{\epsilon}  \vee (\log(n)+\log(N)+c_1 p)\right)^2.
 \end{equation*}{}
 The bound from \cite{park2019non} implies the condition
 \begin{equation*}
     N \geq C \frac{\erank(\mathbf{\Sigma})}{p^{1/2}}\log(2n)\left( \frac{n}{p^{3/2}\epsilon}  \vee (\log(n)+\log(N)+c_1 p)\right)^2,
 \end{equation*}
which requires a significant larger number of samples. This can be observed by the additional $n^2$ factor, and an asymptotic dependency on the missing data in the order of $1/p^{7/2}$.
 \cite{park2019non} points out the necessity for an improved rate of convergence, which is partially addressed by Theorem \ref{th_expected_error}, however since our result is an expectation bound, more research is needed to obtain concentration inequalities. See Table \ref{tab_comparison} for a summary of the sample complexity  comparison from  this section. 
%
%
%
\subsection{\NEW{MCAR case with unknown distribution}}
The main disadvantage of (\ref{eq_est_cov_y_mu}) and Theorem \ref{th_expected_error} is that  in practice, the distribution of $\boldsymbol\delta$ is unknown. Here we present a our bound for (\ref{eq_est_cov_y_mu_mcar2}), which depends on $\widehat{ \mathbf{P}}$ instead of $\mathbf{ P}$.
\begin{theorem}\label{th_expected_error_MAR}
	Suppose the random vector $\mathbf{x}$  satisfies Assumption \ref{assum_trace},  with  $N \geq e^2$ and $n \geq 3$. Also assume that for all $i \in [n]$ there is a $k \in [N]$ so that $\delta_i^{(k)}=1$. Conditioned on the missing data pattern $\mathbf{\Delta}$,  the estimator from (\ref{eq_est_cov_y_mu_mcar2}) satisfies
	\begin{align}\label{spectral_bound_mcar2}
	&\mathbb{E}\left[ \Vert\widehat{\mathbf{\Sigma}} - \mathbf{\Sigma}_{\M}\Vert^2 \mid \mathbf{\Delta } \right]^{1/2} \leq C \Vert \mathbf{\Sigma}\Vert  \\
	& \left( \sqrt{\frac{  8e\log(n)  {\srank}_{\min}(\mathbf{\Sigma},\widehat{\mathbf{P}}) }{N}   } + \frac{4e^2\log(n)\log(N) {\srank}_{2}(\mathbf{\Sigma},\widehat{\mathbf{P}})}{N} \right),\nonumber
	\end{align}
	where $C$ is a  universal numerical constant, $\mathbf{\Sigma}_{\M} = \E[\widehat{\mathbf{\Sigma}}]$, and
	\begin{align*}
	    {\srank}_{\min}(\mathbf{\Sigma},\widehat{\mathbf{P}}) &= \frac{1}{\Vert \mathbf{\Sigma}\Vert} \max_{j}\sum_{i: (i,j)\in \M} \frac{\Sigma_{ii}}{\widehat{p}_{ij}} \\
	    {\srank}_{2}(\mathbf{\Sigma},\widehat{\mathbf{P}}) &= \frac{1}{\Vert \mathbf{\Sigma}\Vert} \sqrt{\sum_{(i,j) \in \M} \frac{\Sigma_{ii}\Sigma_{jj}}{\widehat{p}^2_{ij}}}.
	\end{align*}
\end{theorem}
 In general, (\ref{spectral_bound_mcar2}) corresponds to an  inequality between random variables that are functions of $\mathbf{\Delta}$, however Theorem \ref{th_expected_error_MAR} is also valid if $\mathbf{\Delta}$ is a deterministic matrix.  
The parameters from Theorem \ref{th_expected_error_MAR} are scaled effective ranks, since they obey 
\begin{align}\label{eq_bound1_scaledEffrank_mcar2}
\erank( \mathbf{\Sigma})\leq {\srank}_{\min}(\mathbf{\Sigma},\widehat{\mathbf{P}}) \leq \frac{\erank( \mathbf{\Sigma})}{\widehat{p}_{\min}}, \\
\erank( \mathbf{\Sigma})\leq {\srank}_{2}(\mathbf{\Sigma},\widehat{\mathbf{P}}) \leq \frac{\erank( \mathbf{\Sigma})}{\widehat{p}_{\min}},\label{eq_bound2_scaledEffrank_mcar2}
\end{align} 
almost surely, where $\widehat{p}_{\min} = \min_{(i,j) \in \M} \widehat{p}_{ij}$. 
\paragraph{Comparison with Theorem \ref{th_expected_error}}
 The proposed estimators, (\ref{eq_est_cov_y_mu}) and (\ref{eq_est_cov_y_mu_mcar2}),   can correct   for the bias introduced to the sample covariance matrix due to missing observations, by    normalizing its entries  with the matrices $\mathbf{\Gamma}$ and $\widehat{\mathbf{\Gamma}}$, respectively.  The bound from Theorem \ref{th_expected_error}  is characterized by  scaled effective rank parameters that depend on the observation probability matrix $\mathbf{P}$. In  Theorem \ref{th_expected_error_MAR} we obtained the same bound, however now  the  scaled effective rank parameters have replaced the observation probability matrix with its estimate $\widehat{\mathbf{P}}$. Theorems \ref{th_expected_error} and \ref{th_expected_error_MAR} suggest that, asymptotically (large $N$), both estimators have the same performance, since $\widehat{\mathbf{P}}$ converges to $\mathbf{P}$.
As discussed in Section \ref{sec:miss_estimator}, for $\widehat{\mathbf{\Sigma}}$ to be unbiased, we require that all $\widehat{p}_{ij}>0$, meaning that all pairs of variables must be observed together at least once. Its counterpart in Theorem \ref{th_expected_error} is the condition that all $p_{ij}>0$. 
\paragraph{Comparison with \cite{cai2016minimax}}
   Let us denote the estimator (\ref{eq_estimator_cai_kolar}) by $\tilde{\mathbf{\Sigma}}$, Lemma 2.1 in \cite{cai2016minimax} implies that  with probability at least $1-e^{-\nu}$, $\tilde{\mathbf{\Sigma}}$ obeys
\begin{align*}
\Vert \tilde{\mathbf{ \Sigma}}- \mathbf{ \Sigma}\Vert \leq C \Vert \mathbf{ \Sigma}\Vert  
\left( \sqrt{\frac{c_1 + c_2 n + \nu  }{N \hat{p}_{\min }}} + \frac{c_1 + c_2 n + \nu  }{N \hat{p}_{\min }} \right).
\end{align*}
Essentially, the dominating parameter in the above bound is the quantity $ {n}/{\hat{p}_{\min}}$.  First, this bound depends on the least observed entry (through $\hat{p}_{\min}$), and second, it depends on  $n$ instead of the effective rank, which always satisfies $\erank(\mathbf{\Sigma}) \leq n$.  Therefore, this bound is  
 always larger than the right hand side of (\ref{eq_bound1_scaledEffrank_mcar2}) and (\ref{eq_bound2_scaledEffrank_mcar2}).
\paragraph{A closer look at (\ref{eq_est_cov_y_mu}) and (\ref{eq_est_cov_y_mu_mcar2})}  for (\ref{eq_est_cov_y_mu}), the estimation error of the $ij$ entry is
\begin{equation}\label{eq_entrywise_error_mcar1}
    \E[\vert \widehat{\mathbf{\Sigma}}_{ij} - \mathbf{\Sigma}_{ij}\vert^2] = \frac{\var[(x_i x_j)^2] + (1-p_{ij}) \mathbf{\Sigma}^2_{ij}}{p_{ij}N}.
\end{equation}
A similar calculation produces an entry-wise error for (\ref{eq_est_cov_y_mu_mcar2}) equal to
\begin{equation}\label{eq_entrywise_error_mcar2}
    \E[\vert \widehat{\mathbf{\Sigma}}_{ij} - \mathbf{\Sigma}_{ij}\vert^2] = \E\left[ \frac{1}{\hat{p}_{ij}} \right] \frac{\var[(x_i x_j)^2]}{N}.
\end{equation}
For complete data $p_{ij} = \hat{p_{ij}}=1$,  both estimators have the same error. When there is missing data, if  $N$ is large, $\E\left[ {1}/{\hat{p}_{ij}} \right]$ converges to $1/p_{ij}$, and thus  (\ref{eq_est_cov_y_mu}) will always be outperformed by (\ref{eq_est_cov_y_mu_mcar2}). Note that since $\E[x_ix_j]^2 = \mathbf{\Sigma}_{ij}^2$, (\ref{eq_entrywise_error_mcar2}) is composed of a variance term, while (\ref{eq_entrywise_error_mcar1}) has a variance and an additional bias. This occurs because the $ij$ entries of both estimators are a sum of $\hat{p}_{ij}N$ terms, however since for  (\ref{eq_est_cov_y_mu}), the number of terms does not coincide with the dividing factor $p_{ij}N$, each term is not an exact average, unlike the case of (\ref{eq_est_cov_y_mu_mcar2}).
\subsection{CMCAR case with known distribution}
Now we present an error bound for   (\ref{eq_cov_est_martinga}). 
\begin{theorem}\label{th_expected_error_martin} Suppose  $\mathbf{P}^{(t+1)}=f_t(\widehat{\mathbf{\Sigma}}_t)$ for all $t \geq 1$ where $f_t(~\cdot~)$ are bounded, in the sense that there is a constant $\alpha $ so that for all $i$ and for all $t$, $p_{ij}^{(t)} \geq \alpha>0$ almost surely. Under the same assumptions of Theorem \ref{th_expected_error},    (\ref{eq_cov_est_martinga}) obeys
\begin{align*}
\mathbb{E}[\Vert \widehat{\mathbf{\Sigma}}_T - \mathbf{\Sigma} \Vert^2]^{1/2} &\leq \tilde{C} \Vert \mathbf{\Sigma}\Vert \sqrt{ \log(n+1)}\\
&\times\left[ \sqrt{ \frac{ 8e\log(n) }{N} \varepsilon_1} + \frac{4e^2\log(n)L}{N}\varepsilon_2\right],
\end{align*}
where $\tilde{C}$ is an universal constant,   $ \varepsilon_1, \varepsilon_2$ are  given by
\NEW{
\begin{align*}
\varepsilon_1(\mathbf{\Sigma},\lbrace \mathbf{P}^{(t)},N_t \rbrace_{t=1}^T)&=\E\left[\frac{1}{N}\sum_{t=1}^T {{\srank}_{\min}}(\mathbf{\Sigma},\mathbf{P}^{(t)}) N_t\right] \\
\varepsilon_2^2(\mathbf{\Sigma},\lbrace \mathbf{P}^{(t)},N_t \rbrace_{t=1}^T)&=\E\left[\frac{1}{L^2}\sum_{t=1}^T {{\srank}_2}^2(\mathbf{\Sigma},\mathbf{P}^{(t)})\log^2(N_t)\right],
\end{align*}
}
and  $N = \sum_{t=1}^T N_t$, and $L^2 = \sum_{t=1}^T \log^2(N_t)$.
\end{theorem}
Note that when $T=1$, the bounds from Theorems \ref{th_expected_error} and \ref{th_expected_error_martin} only differ on the $\sqrt{\log(n+1)}$ term and a constant. 
The estimation error depends on the behavior of the quantities $\varepsilon_i$, which depends on the observations,    the missing data mechanism and $\lbrace  f_t \rbrace_{t=1}^T$. In particular, when $f_t$ depends only on the index $t$, and not on the data, we have the case of a time varying MCAR  mechanism.
The assumptions of Theorem \ref{th_expected_error_martin} guarantee that $\varepsilon_i$ are  scaled effective ranks, thus
\begin{equation*}
\varepsilon_i \leq  \frac{\erank(\mathbf{\Sigma})}{\alpha}, 
\end{equation*}
and $\widehat{\mathbf{\Sigma}}_T \rightarrow \mathbf{\Sigma}$ in probability.
To  establish that this bound is also consistent  with Theorems \ref{th_expected_error} and \ref{th_expected_error_MAR}  for $T>1$ we will analyze the limiting behavior of the expectations in the upper bound.
\begin{proposition}\label{prop_convergence_objective}
	Under the same assumptions of Theorem \ref{th_expected_error_martin},  if  additionally $f_t=f$ is  continuous,   then
	\begin{align}
	&\mathbf{P}^{(t)} \rightarrow \mathbf{P}^{\infty}, \\
	&{\srank}_2(\mathbf{\Sigma}, \mathbf{P}^{(t)}) \rightarrow {\srank}_2(\mathbf{\Sigma}, \mathbf{P}^{\infty}), \\
	&{\srank}_{\min}(\mathbf{\Sigma}, \mathbf{P}^{(t)}) \rightarrow {\srank}_{\min}(\mathbf{\Sigma}, \mathbf{P}^{\infty}),
	\end{align}
	in probability, where $\mathbf{P}^{\infty} = f(\mathbf{\Sigma})$.
	Also we have convergence of the expectation, that is,
	\begin{align}
	\E[{\srank}^2_2(\mathbf{\Sigma}, \mathbf{P}^{(t)})] \rightarrow {\srank}_2^2(\mathbf{\Sigma}, \mathbf{P}^{\infty}),\\
	\E[{\srank}_{\min}(\mathbf{\Sigma}, \mathbf{P}^{(t)})] \rightarrow {\srank}_{\min}(\mathbf{\Sigma}, \mathbf{P}^{\infty}).
	\end{align}
	\begin{IEEEproof}
		Since $f$ is bounded and  a continuous function of $\widehat{\mathbf{\Sigma}}_t$, so  are ${\srank}_s \circ f$, and ${\srank}_{\min} \circ f$. All limits are implied by convergence in probability of  $\widehat{\mathbf{\Sigma}}_t$ to its expected value, along with  \cite[Th. 2.3.4]{durrett2010probability}.
	\end{IEEEproof}
\end{proposition}
Proposition \ref{prop_convergence_objective}, and  Theorems \ref{th_expected_error} and \ref{th_expected_error_martin}  imply that asymptotically, our  estimators in (\ref{eq_est_cov_y_mu}) and  (\ref{eq_cov_est_martinga}) converge at the same rates. This is a direct consequence of the following limits, which are implied by Proposition \ref{prop_convergence_objective}
\begin{align}
\varepsilon_1 &\rightarrow {\srank}_{\min}(\mathbf{\Sigma}, \mathbf{P}^{\infty}), \textnormal{ and }\\
\varepsilon_2 &\rightarrow {\srank}_{2}(\mathbf{\Sigma}, \mathbf{P}^{\infty}).
\end{align}
We will finish the analysis of this estimator by establishing a lower bound on the sample complexity. 
\begin{proposition}\label{prop_sample_complex_martin}
	Suppose that  $e < N_t$ for all $t \in [T]$, and  the following conditions hold
	\begin{align}
	&\E\left[ \sum_{t=1}^T {\srank}_{\min}(\mathbf{\Sigma}, \mathbf{P}^{(t)}) N_t \right]\nonumber  \\
	&\geq \E\left[ \sum_{t=1}^T2e^3\log(n)\log^2(N_t){\srank}^2_{2}(\mathbf{\Sigma}, \mathbf{P}^{(t)})  \right], \label{eq_condition_Nt}\\
	N &\geq 32e\tilde{C}^2\log(n)\log(n+1)\varepsilon_1 \frac{1}{\epsilon^2} \label{eq_condition_sample_comp_martin}
	\end{align}
	for some $T >1$ and some $\epsilon>0$, then
	\begin{equation*}
	\mathbb{E}\left[ \Vert\widehat{\mathbf{\Sigma}}_T - \mathbf{\Sigma}\Vert^2 \right]^{1/2}  \leq  \epsilon \Vert\mathbf{\Sigma}\Vert.
	\end{equation*}
	\begin{IEEEproof}
		The first inequality implies
		\begin{equation*}
		\sqrt{ \frac{ 8e\log(n) }{N} \varepsilon_1} \geq  \frac{2e^3\log(n)L}{N}\varepsilon_2.
		\end{equation*}
		Combining this with Theorem \ref{th_expected_error_martin} leads to
		\begin{equation*}
		\mathbb{E}[\Vert \widehat{\mathbf{\Sigma}}_T - \mathbf{\Sigma} \Vert^2]^{1/2} \leq 2\tilde{C} \Vert \mathbf{\Sigma}\Vert 
		\left[ \sqrt{ \frac{ 8e\log(n)\log(n+1) }{N} \varepsilon_1} \right].
		\end{equation*}
		Using the second assumption leads to the desired bound.
	\end{IEEEproof}
\end{proposition}
If we ignore   logarithmic and constant factors the condition (\ref{eq_condition_sample_comp_martin}) would be in the same order of magnitude as the sample complexity we established  for the MCAR estimators. Note that in this case there is an additional $\log(n+1)$ factor.
The inequality (\ref{eq_condition_Nt}) is a condition on the sequence $\lbrace N_t \rbrace$, and although it looks convoluted, it is easily satisfied if $T$ is large enough, and 
\begin{equation}
N_t/\log^2(N_t) = \Omega(\log(n) {\srank}_{2}^2(\mathbf{\Sigma}, \mathbf{P}^{\infty })/{\srank}_{\min}(\mathbf{\Sigma}, \mathbf{P}^{\infty })).
\end{equation}
This can be achieved for example if  ${N_t}$ is an increasing sequence with $\lim_{t \rightarrow \infty} N_t = \infty$. 

\section{Experiments}
 \label{sec:miss_exp_cov}
\subsection{Covariance estimation under MCAR observations}
\subsubsection{Uniform MCAR observations}
\label{exp_mcar_uniform}
We first compare the estimators from (\ref{eq_est_cov_y_mu}) and (\ref{eq_est_cov_y_mu_mcar2}),  designed to handle MCAR observations with known and unknown missing data distribution respectively.  The experimental setup is as follows.
\begin{itemize}
	\item Data follows a Gaussian distribution with mean zero,  and  population covariance matrix $\mathbf{\Sigma}$. The eigenvalues are $\rho^i$ for $i=1,\cdots,n$, with $0<\rho<1$ chosen so that $\erank({\mathbf{\Sigma}})=4$. The eigenvectors are generated by applying the Gram-Schmidt method to a $n \times n$ random matrix with $\mathcal{N}(0,1)$ independent entries.
	\item The dimension is set to $n=50$. The MCAR missing data mechanism is uniform with independend entries and  observation probabilities $p \in \lbrace  0.4, 0.6 , 0.8\rbrace $.
	\item We generate complete and incomplete data. The number of samples     $N$ are chosen to be  $100$ equidistant  values in the log scale between $15$ and  $50n$.
	\item For each pair $(N,p)$ we compute estimation error and average over $100$ trials.
\end{itemize}

\begin{figure}[htbp]
	\centering
	\includegraphics[width=0.48\textwidth]{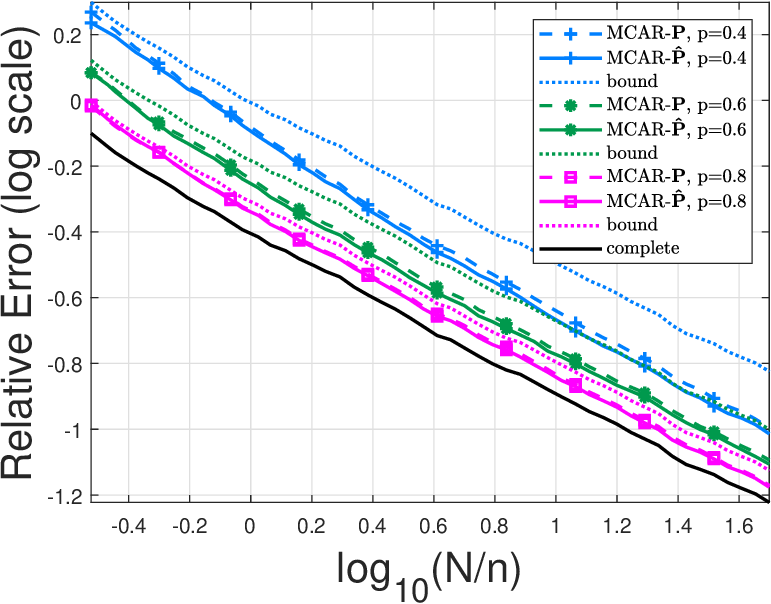}
	\caption{Comparison of proposed estimators under MCAR missing observations with uniform observation probabilities. The plot shows relative estimation error in operator norm versus number of i.i.d. realizations, both in logarithmic scale. Solid and dashed lines denote estimators (\ref{eq_est_cov_y_mu}) (MCAR-$\widehat{\mathbf{P}}$) and (\ref{eq_est_cov_y_mu_mcar2}) (MCAR-${\mathbf{P}}$) respectively. Dotted lines are an estimate based on the error bounds. }
	\label{fig_mcar_vs_mar}
\end{figure}

The averaged estimation errors are plotted in Figure \ref{fig_mcar_vs_mar}.  The black line in the figure corresponds to the estimation error of the sample covariance matrix, and serves as our baseline performance.
Theorem \ref{th_expected_error} indicates that for $N$ large enough, the MCAR estimator with missing data has error $\mathcal{O}((p \sqrt{N})^{-1})$, thus  to achieve the same error as the sample covariance matrix, it needs more samples by a factor of $p^{-1}$. We verify this experimentally by multiplying the estimation error of the sample covariance matrix by $p^{-1}$. Those curves are  shown in dotted lines, and are always above the estimation error of both  estimators, almost parallel and tighter when $N/n<1$. Additionally,  when $\mathbf{P}$ is unknown,  (\ref{eq_est_cov_y_mu_mcar2}) always outperforms  (\ref{eq_est_cov_y_mu}), even though the latter has access to the actual missing data distribution through $\mathbf{P}$. This might occur because  (\ref{eq_est_cov_y_mu_mcar2}) computes exact averages (using $\widehat{\mathbf{P}}$), while (\ref{eq_est_cov_y_mu}) employs an   estimate of the number of observations (a function of $\mathbf{P}$). 
\subsubsection{Non uniform MCAR observations}
\begin{figure}[htbp]
	\centering
	\begin{subfigure}[b]{0.34\textwidth}
		\centering
		\includegraphics[width=1\textwidth]{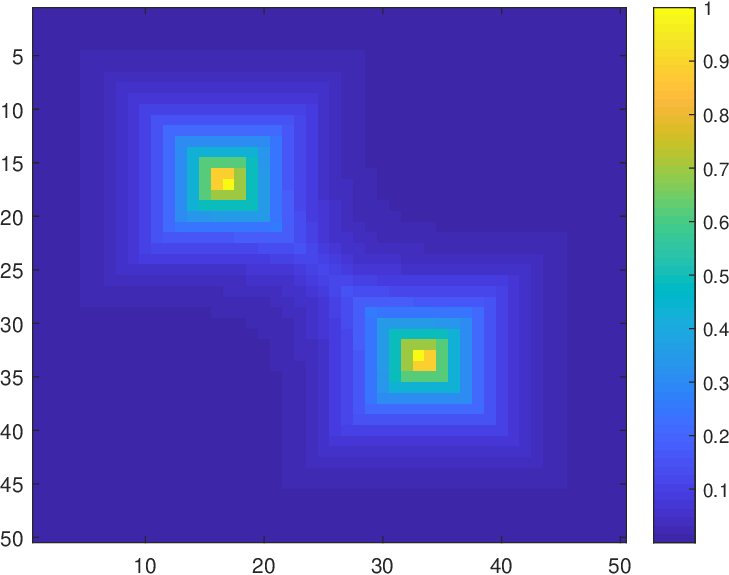}
		\caption{covariance matrix}
		\label{fig_cov_skewed:cov}
	\end{subfigure}
	\begin{subfigure}[b]{0.34\textwidth}
		\centering
		\includegraphics[width=1\textwidth]{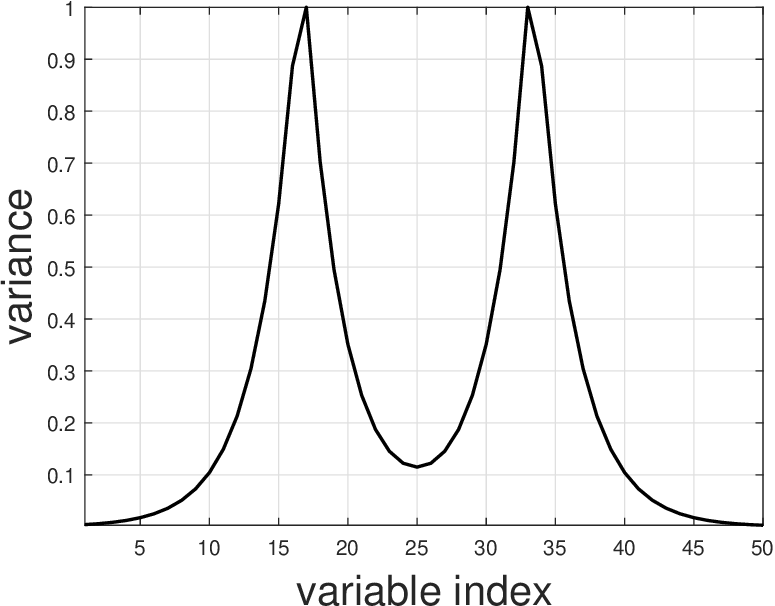}
		\caption{variance}
		\label{fig_cov_skewed:plotvars}
	\end{subfigure}
	\caption{Approximately sparse covariance matrix with low effective rank. }
	\label{fig_cov_skewed}
\end{figure}

\begin{figure}[htbp]
	\centering
	\includegraphics[width=0.48\textwidth]{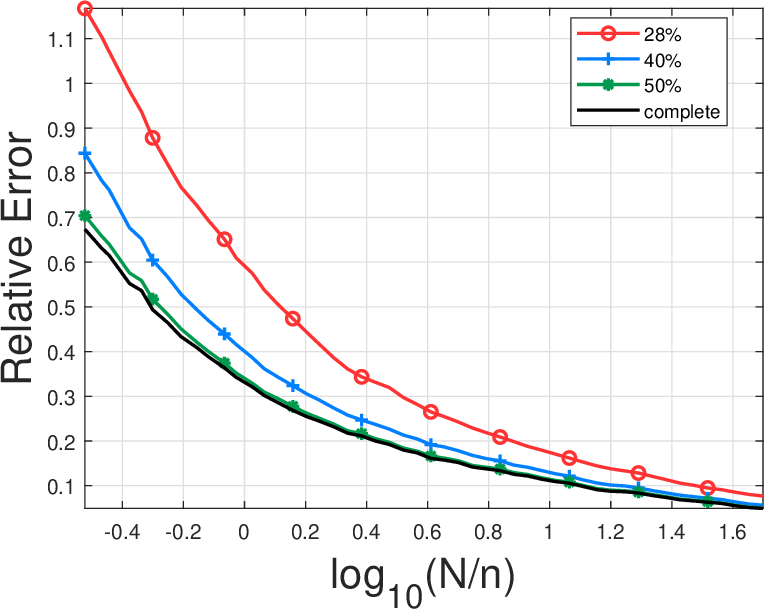}
	\caption{Evaluation of the  estimators under MCAR missing observations with non uniform observation probabilities.}
	\label{fig_mar_favorable_p}
\end{figure}

In the  next experiment we show  that when the missing observations have a favorable distribution, not much information is lost due to missing data. We consider a non-singular population covariance matrix $\mathbf{ \Sigma}$ of dimension $n=50$, and effective rank $\erank = 2.7472$ (see Figure \ref{fig_cov_skewed}). 
The missing data distribution is MCAR with independent entries, and observation probabilities given by $\mathbf{p}= \max(\rho \sqrt{\diag(\mathbf{\Sigma})}, 1)$. In this setting the variable $x_i$ is observed with probability proportional to its standard deviation $\sqrt{\Sigma_{ii}}$. The parameter $\rho$ is set so that the expected number of observed variables are $\lbrace 28\%, 40\%, 50\% \rbrace $. The performance of (\ref{eq_est_cov_y_mu_mcar2}) is depicted  in Figure \ref{fig_mar_favorable_p}. It can be seen that when roughly $50\%$ or more  variables are observed in average, the estimator has almost the same performance as the sample covariance matrix that has access to complete observations. This can be attributed to the fact that variables with lower variance contribute less to the estimation error, and thus a higher missing data rate can be tolerated for those  entries. 
%
\begin{figure}[ht]
	\centering
	\includegraphics[width=0.48\textwidth]{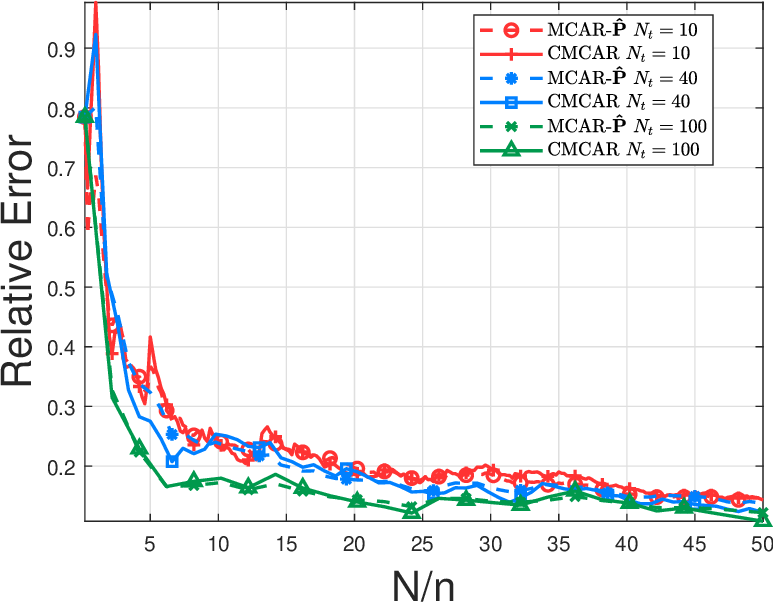}
	\caption{Comparison of CMCAR and MCAR estimators under  time varying  missing observations with uniform observation probabilities. The plotted curve shows relative estimation error (in log scale) in operator norm versus number of i.i.d. realizations. $N_t$ is the block size, i.e., the number of realizations until the missing data distribution changes.  }
	\label{fig_cmcar_vs_mar}
\end{figure}
\subsection{Estimation under CMCAR observations}
This section compares the estimators from (\ref{eq_est_cov_y_mu_mcar2}) and (\ref{eq_cov_est_martinga}), designed for  MCAR and CMCAR observations respectively. We consider a simple CMCAR mechanism, where  
\begin{equation*}
    \Prob(\delta_i^{(t,k)}=1) = (1-\alpha_t) + \alpha_t \frac{m}{n},
\end{equation*}{}
for  fixed $m$ and   $    \alpha_t = 0.9\left( 1-{1}/{\sqrt{t}}\right).$
When $t=1$, there are complete observations, and the amount of missing data increases with $t$, to reach an uniform observation probability of $0.9m/n$. We consider the same synthetic dataset created in Section \ref{exp_mcar_uniform}, and set $m=0.1n$. We fix $N_1=10$ and the remaining block sizes  $N_t$ for $t\geq 2$ are set to a constant.  Figure \ref{fig_cmcar_vs_mar} depicts the average error in operator norm of both estimators. Note that the  CMCAR estimator (\ref{eq_cov_est_martinga}) knows the missing data distribution, while the MCAR estimator (\ref{eq_est_cov_y_mu_mcar2}) does not. 
When the missing data distribution changes faster (smaller $N_t$), both estimators require more samples to be accurate. Moreover, when the block size increases, and the missing data distribution changes more slowly, the estimation error decays faster, and both estimators behave similarly. 
%
\section{Proofs}
 \label{sec:miss_proofs}
\subsection{Bounding tools}
\paragraph{Sub-exponential and sub-Gaussian properties  \cite{vershynin2016high}}
\begin{itemize}
	\item If $z$ is a sub-Gaussian random variable, then
	\begin{equation*}
	\E [\vert z\vert ^r] ^{1/r} \leq c \Vert z \Vert_{\Psi_2} \sqrt{r}, \textnormal{ for } r\geq 1.
	\end{equation*}
	\item If $z$ and $w$ are sub-Gaussian random variables, then $z^2$ and $zw$ are sub-exponential with norms satisfying $\Vert z^2 \Vert_{\Psi_1} = \Vert z \Vert_{\Psi_2}^2$, and $\Vert zw \Vert_{\Psi_1} \leq \Vert z \Vert_{\Psi_2} \Vert w \Vert_{\Psi_2}$.
	\item Let $y_1 = \delta_1 x_1$ and $y_2 = \delta_2 x_2$, where  $\delta_1, \delta_2$ are Bernoulli with  parameters $p_1$ and $p_2$, and  $x_1, x_2$ are sub-Gaussian random variables. If $\delta_1, \delta_2$ are independent of  $x_1$ and $x_2$, then
	\begin{enumerate}
		\item[1)]    $y_i$ is sub-Gaussian and $\Vert y_i \Vert_{\Psi_2} 	\leq \Vert x_i \Vert_{\Psi_2}$.
		\item[2)]  $y_1^2$, $y_2^2$, and $y_1 y_2$ are sub-exponential with norms satisfying  $\Vert y_i y_j \Vert_{\Psi_1} \leq \Vert x_i x_j \Vert_{\Psi_1}$ for $i,j =1,2$.
	\end{enumerate}
\end{itemize}
\paragraph{Useful  inequalities}
\begin{itemize}
	\item Jensen's inequality. Let $z$ be a real random variable, with finite $r$-th moment, then
	\begin{equation*}
	\E[\vert z\vert]^{r} \leq \E [\vert z \vert ^r], \textnormal{ for } r\geq 1.
	\end{equation*}
	\item Application of Jensen's inequality.  Let $z_1,\cdots, z_N$ be  non negative random variables with finite $r$-th moments,  then 
	\begin{equation}\label{eq_expected_max}
	\E[\max_{k} z_k] \leq N^{1/r} \max_k\E[z_k^r]^{1/r},  \textnormal{ for } r\geq 1.
	\end{equation}
	\item Minkowski's inequality (triangular inequality). Given  $z_1,z_2$ with finite $r$-th moments, then for any $r\geq 1$
	\begin{equation*}
	\E[\vert z_1+z_2 \vert ^r]^{1/r} \leq 		\E[\vert z_1 \vert ^r]^{1/r} + 	\E[\vert z_2 \vert ^r]^{1/r}.
	\end{equation*}
	\item For any $a,b \in \mathbb{R}$,  $(a+b)^2 \leq 2(a^2 + b^2) $.
\end{itemize} 
\paragraph{Hadamard products}
\begin{itemize}
	\item We will need the following identify
	\begin{equation}\label{eq_squared_hadamard}
	\left(\mathbf{ M \odot  x x^{\top}}\right)^2 = \sum_{i=1}^n x_i^2 \mathbf{x x^{\top}} \odot \mathbf{M}_i \mathbf{M}_i^{\top },
	\end{equation}
	where $\mathbf{M}_i$ is the $i$-th column of $\mathbf{M}$.
	\item If  $\mathbf{A}_1 \prec \mathbf{A}_2$, and $\mathbf{B}$ is  positive semi-definite, then
	\begin{equation}\label{eq_schur_theorem}
	\mathbf{A}_1 \odot \mathbf{B} \prec 	\mathbf{A}_2 \odot \mathbf{B}
	\end{equation}
	where  $\mathbf{A} \prec \mathbf{B}$ stands for $0 \prec \mathbf{B-A}$. In particular we have that
	\begin{equation}
	\mathbf{A}_1 \odot \mathbf{B} \prec 	\Vert \mathbf{A}_1\Vert \mathbf{I} \odot \mathbf{B} .
	\end{equation}
\end{itemize}
Some of our  proofs rely on an expectation bound from \cite{chen2012masked}, originally devised  to certify the  performance of  the masked covariance estimator. By  following Theorem 3.2 and the beginning of the proof of Theorem 3.1 in  \cite{chen2012masked} we state their result in the form most useful for us.
\begin{theorem}[\cite{chen2012masked}] \label{theorem_chen_masked}
	Let $\lbrace \mathbf{Z}_k\rbrace$ be a finite sequence of independent, symmetric, self adjoint random matrices of size $n \times n$, and $n \geq 3$, then
	\begin{align*}
	&\mathbb{E}\left[ \left\Vert \sum_k \mathbf{Z}_k - \E[\mathbf{Z}_k]  \right\Vert^2 \right]^{1/2} \leq \\ &\sqrt{\alpha} \left\Vert  \sum_k \E[\mathbf{Z}_k^2] \right\Vert^{1/2} + 
	\alpha \E\left[ \max_k \Vert \mathbf{Z}_k \Vert^2 \right]^{1/2},
	\end{align*}
	where $\alpha  = 8e\log(n)$.
\end{theorem}
\subsection{Proof of Theorem \ref{th_expected_error} }
Applying Theorem \ref{theorem_chen_masked} to (\ref{eq_est_cov_y_mu})  produces
\begin{align*}
&\mathbb{E}\left[ \Vert\widehat{\mathbf{\Sigma}} - \mathbf{\Sigma}\Vert^2 \right]^{1/2} \leq \sqrt{\frac{8e\log(n)\Vert \mathbb{E}\left[ (\mathbf{\Gamma} \odot \mathbf{y} \mathbf{y}^{\top} )^2 \right]\Vert}{N}}  \\
&+\frac{8e\log(n)\mathbb{E}\left[ \max_{k}\Vert\mathbf{\Gamma} \odot \mathbf{y}^{(k)} {\mathbf{y}^{(k)}}^{\top} \Vert^2 \right]^{1/2}}{N}.
\end{align*}
The bounds  from \cite{chen2012masked} for $\Vert \mathbb{E}\left[ (\mathbf{\Gamma} \odot \mathbf{y} \mathbf{y}^{\top} )^2 \right]\Vert$ and $\mathbb{E}\left[ \max_{k}\Vert\mathbf{\Gamma} \odot \mathbf{y}^{(k)} {\mathbf{y}^{(k)}}^{\top} \Vert^2 \right]$  are fairly general. They   decouple $\mathbf{y}$ and $\mathbf{y}^{(k)}$ from $\mathbf{\Gamma}$, thus they  can be greatly improved if we consider their  relation   for this  problem. 
\paragraph{Bounding the matrix variance (first term)}
Using  (\ref{eq_squared_hadamard}) and applying expectation we obtain
\begin{equation}\label{eq_squared_hadamard_mcar}
	\E\left[\left(\mathbf{ \Gamma \odot  y y^{\top}}\right)^2 \right]= \sum_{i=1}^n \E[ y_i^2 \mathbf{y y^{\top}}] \odot \mathbf{\Gamma}_i \mathbf{\Gamma}_i^{\top }.
	\end{equation}
Using independence, each term in the expected value above can be decomposed as follows
\begin{equation*}
    \E[ y_i^2 \mathbf{y y^{\top}}] = \E[ x_i^2 \delta_i \mathbf{x x^{\top}} \odot \boldsymbol\delta \boldsymbol\delta^{\top}] =    \E[ x_i^2  \mathbf{x x^{\top}} ]  \odot \E[ \delta_i \boldsymbol\delta \boldsymbol\delta^{\top}] .
\end{equation*}
Replacing back in (\ref{eq_squared_hadamard_mcar}) and applying  (\ref{eq_schur_theorem}) to the Hadamard product of positive semi-definite matrices we have
\begin{align*}
	\E\left[\left(\mathbf{ \Gamma \odot  y y^{\top}}\right)^2 \right]&= \sum_{i=1}^n  \E[ x_i^2  \mathbf{x x^{\top}} ]  \odot \E[ \delta_i \boldsymbol\delta \boldsymbol\delta^{\top}] \odot \mathbf{\Gamma}_i \mathbf{\Gamma}_i^{\top } \\
	&\prec \sum_{i=1}^n  \left\Vert \E[ x_i^2  \mathbf{x x^{\top}} ]\right\Vert \mathbf{I }  \odot \E[ \delta_i \boldsymbol\delta \boldsymbol\delta^{\top}] \odot \mathbf{\Gamma}_i \mathbf{\Gamma}_i^{\top }.
	\end{align*}
Since  the right side is a  diagonal matrix, its operator norm is its largest entry, and   we can bound the matrix variance as
\begin{align*}
&\left\Vert 	\E\left[\left(\mathbf{ \Gamma \odot  y y^{\top}}\right)^2 \right] \right\Vert \leq \\
&\max_{j \in [n]}\sum_{i=1}^n  \left\Vert \E[ x_i^2  \mathbf{x x^{\top}} ]\right\Vert \left( \E[ \delta_i \boldsymbol\delta \boldsymbol\delta^{\top}]  \mathbf{\Gamma}_i \mathbf{\Gamma}_i^{\top } \right)_{jj}.
\end{align*}
The following holds
\begin{equation*}
\left( \E[ \delta_i \boldsymbol\delta \boldsymbol\delta^{\top}]  \mathbf{\Gamma}_i \mathbf{\Gamma}_i^{\top } \right)_{jj} = \frac{1}{p_{ij}}.
\end{equation*}
The desired bound is obtained after applying Lemma \ref{lemma_auxiliary_1} to bound $\left\Vert \E[ x_i^2  \mathbf{x x^{\top}} ]\right\Vert $ (at the end of this section), thus producing
\begin{equation*}
\left\Vert 	\E\left[\left(\mathbf{ \Gamma \odot  y y^{\top}}\right)^2 \right] \right\Vert \leq C^2 \Vert \mathbf{\Sigma}\Vert\max_{j \in [n]}	  \sum_{i=1}^n \frac{ \Sigma_{ii}}{p_{ij}},
\end{equation*}
where  $C=4c^2/c_0$.
When $\boldsymbol\delta$ has  independent coordinates, we have that $p_{ij} = p_{ii}p_{jj}$ when $i \neq j$, therefore we may use the bound
\begin{equation*}
\Vert \mathbb{E}\left[ (\mathbf{\Gamma} \odot \mathbf{y} \mathbf{y}^{\top} )^2 \right]  \Vert \leq C^2\frac{\Vert \mathbf{\Sigma}\Vert}{p_{\min}}\left(  \sum_{k=1}^n \frac{\Sigma_{kk}}{p_k} \right),
\end{equation*}
where $p_{\min} = \min_{i}p_{ii}$.
\paragraph{Bounding the maximum spectral norm (second term).}
For any $r \geq 1$,   applying (\ref{eq_expected_max}) to $\Vert\mathbf{\Gamma} \odot \mathbf{y}^{(k)} {\mathbf{y}^{(k)}}^{\top} \Vert^2$ results in
\begin{equation*}
\mathbb{E}\left[ \max_{k}\Vert\mathbf{\Gamma} \odot \mathbf{y}^{(k)} {\mathbf{y}^{(k)}}^{\top} \Vert^2  \right]\leq N^{1/r} \mathbb{E}\left[ \Vert\mathbf{\Gamma} \odot \mathbf{y} {\mathbf{y}}^{\top} \Vert^{2r } \right]^{1/r}.
\end{equation*}
The dependency on $k$ is removed because $\mathbf{y}^{(k)}$ are i.i.d. copies of $\mathbf{y}$.
Now we only need to bound $\mathbb{E}\left[ \Vert\mathbf{\Gamma} \odot \mathbf{y} {\mathbf{y}}^{\top} \Vert^{2r } \right]$.Since the operator norm is bounded by the Frobenius norm, and  $y_i^2 \leq x_i^2$,  we have the inequalities
\begin{equation*}
\Vert  \mathbf{\Gamma} \odot \mathbf{y} {\mathbf{y}}^{\top}\Vert^2 
\leq  \Vert  \mathbf{\Gamma} \odot \mathbf{y} {\mathbf{y}}^{\top}\Vert_F^2 
= \sum_{i,j}\gamma_{ij}^2 y_i^2 y_j^2 
\leq \sum_{i,j}\gamma_{ij}^2 x_i^2 x_j^2.
\end{equation*}
Now we apply the triangular and Cauchy-Schwartz inequalities, followed by Assumption \ref{assum_trace}, thus
\begin{align*}
\mathbb{E}\left[ \Vert\mathbf{\Gamma} \odot \mathbf{y} {\mathbf{y}}^{\top} \Vert^{2r } \right]^{1/r} &\leq \sum_{i,j} \gamma_{ij}^2\E\left[ (x_i^2 x_j^2)^r\right]^{1/r} \\
&\leq \sum_{i,j} \gamma_{ij}^2 \E\left[ x_i^{4r}\right]^{1/2r} \E\left[ x_j^{4r}\right]^{1/2r} \\
&\leq (4rc^2/c_0)^2  \sum_{i,j} \gamma_{ij}^2 \Sigma_{ii} \Sigma_{jj}.
\end{align*}
Replacing $\gamma_{ij} = {1}/{p_{ij}}$, leads to
\begin{equation*}
\mathbb{E}\left[ \Vert\mathbf{\Gamma} \odot \mathbf{y} {\mathbf{y}}^{\top} \Vert^{2r } \right]^{1/r} \leq  C^2 r^2  \sum_{i=1}^n\sum_{j=1}^n \frac{\Sigma_{ii}\Sigma_{jj}}{p^2_{ij}}. 
\end{equation*}
The final bound is obtained by minimizing the  function $N^{1/r} r^2$ over $r\geq 1$. This is achieved by setting $r = \log(N)/2$,  since the assumption $\log(N) \geq 2$ guarantees that $r \geq 1$ is satisfied. The desired bound is 
\begin{equation*}
	\mathbb{E}\left[ \max_{k}\Vert\mathbf{\Gamma} \odot \mathbf{y}^{(k)} {\mathbf{y}^{(k)}}^{\top} \Vert^2  \right]\leq   \left(e C\frac{\log(N)}{2}\right)^2  \sum_{i,j=1}^n \frac{\Sigma_{ii}\Sigma_{jj}}{p^2_{ij}}.
\end{equation*}
When $\boldsymbol\delta $ has independent coordinates, we may use the bound
\begin{equation*}
	\mathbb{E}\left[ \max_{k}\Vert\mathbf{\Gamma} \odot \mathbf{y}^{(k)} {\mathbf{y}^{(k)}}^{\top} \Vert^2  \right]\leq  \left(e C\frac{\log(N)}{2} \sum_{i=1}^n \frac{\Sigma_{ii}}{p_{ii}^2} \right)^2.
\end{equation*}
\subsection{Proof of Theorem \ref{th_expected_error_MAR}}
Since  $\E[\widehat{\mathbf{ \Sigma}} \mid \mathbf{\Delta}]=\mathbf{\Sigma}_{\M}$, we can apply Theorem \ref{theorem_chen_masked} to (\ref{eq_est_cov_y_mu_mcar2}), 
\begin{align*}
&\mathbb{E}\left[ \Vert\widehat{\mathbf{\Sigma}} - \mathbf{\Sigma}_{\M}\Vert^2 \mid \mathbf{\Delta} \right]^{1/2} \leq \\ &\frac{1}{N}\sqrt{8e\log(n)\left\Vert \sum_{k=1}^N\mathbb{E}\left[ (\widehat{\mathbf{\Gamma}} \odot \mathbf{y}^{(k)} {\mathbf{y}^{(k)}}^{\top} )^2 \mid \mathbf{\Delta}  \right]\right\Vert}  \\
&+\frac{8e\log(n)\mathbb{E}\left[ \max_{k}\Vert \widehat{\mathbf{\Gamma}} \odot \mathbf{y}^{(k)} {\mathbf{y}^{(k)}}^{\top} \Vert^2 \mid \mathbf{\Delta}  \right]^{1/2}}{N}.
\end{align*}
We follow a   strategy close to the one used in the proof of Theorem \ref{th_expected_error}. Here,   the random vectors $\mathbf{y}^{(k)}$ are conditionally  independent, given the missing data pattern $\mathbf{\Delta}$. In addition, these vectors are  not identically distributed, therefore computations are more cumbersome.
\paragraph{Bounding the matrix variance}
 To simplify computations  we define $\mathbf{M}^{(k)} = \widehat{\mathbf{\Gamma}} \odot \boldsymbol\delta^{(k)} {\boldsymbol\delta^{(k)}}^{\top}$, therefore 
\begin{equation*}
\widehat{\mathbf{\Gamma}} \odot \mathbf{y}^{(k)} {\mathbf{y}^{(k)}}^{\top} = \mathbf{M}^{(k)} \odot \mathbf{x}^{(k)} {\mathbf{x}^{(k)}}^{\top} .
\end{equation*}
\NEW{Since  $\mathbf{M}^{(k)}$ depends only on $\mathbf{\Delta}$, we can decouple the data and the missing data quantities. The following equalities are based on the MCAR property, and the i.i.d.   assumption on $\mathbf{x}^{(k)}$. We obtain the desired conditional expectation after applying  (\ref{eq_squared_hadamard}), resulting in}
\begin{align*}
&\sum_{k=1}^N \mathbb{E}\left[ (\mathbf{M}^{(k)} \odot \mathbf{x}^{(k)} {\mathbf{x}^{(k)}}^{\top} )^2  \mid \mathbf{\Delta} \right] \\
&= \sum_{k=1}^N \sum_{i=1}^n \E[{x_i^{(k)}}^2\mathbf{x}^{(k)} {\mathbf{x}^{(k)}}^{\top}  \mid \mathbf{\Delta}] \odot \mathbf{ M}^{(k)}_i {\mathbf{ M}_i^{(k)}}^{\top} \\
&=\sum_{k=1}^N \sum_{i=1}^n \E[{x_i}^2\mathbf{x} {\mathbf{x}}^{\top} \mid \mathbf{\Delta}] \odot \mathbf{ M}^{(k)}_i {\mathbf{ M}_i^{(k)}}^{\top} \\
&= \sum_{i=1}^n \E[{x_i}^2\mathbf{x} {\mathbf{x}}^{\top} ] \odot \sum_{k=1}^N\mathbf{ M}^{(k)}_i {\mathbf{ M}_i^{(k)}}^{\top}.
\end{align*}
Now we use (\ref{eq_schur_theorem}) to bound the Hadamard product of positive definite matrices, thus  we have
\begin{align*}
&\sum_{k=1}^N \mathbb{E}\left[ (\mathbf{M}^{(k)} \odot \mathbf{x}^{(k)} {\mathbf{x}^{(k)}}^{\top} )^2  \mid \mathbf{\Delta} \right]\\
&\prec \sum_{i=1}^n \Vert \E[{x_i}^2\mathbf{x} {\mathbf{x}}^{\top}]  \Vert \mathbf{I}\odot\sum_{k=1}^N\mathbf{ M}^{(k)}_i {\mathbf{ M}_i^{(k)}}^{\top}.
\end{align*}
The resulting matrix on the right side is diagonal, hence its  operator norm can be computed in closed form, and the desired operator norm can be bounded as
\begin{align*}
&\left\Vert \sum_{k=1}^N \mathbb{E}\left[ (\mathbf{M}^{(k)} \odot \mathbf{x}^{(k)} {\mathbf{x}^{(k)}}^{\top} )^2  \mid \mathbf{\Delta} \right] \right\Vert \\
&\leq \left\Vert  \sum_{i=1}^n \Vert \E[{x_i}^2\mathbf{x} {\mathbf{x}}^{\top}]  \Vert \mathbf{I}\odot\sum_{k=1}^N\mathbf{ M}^{(k)}_i {\mathbf{ M}_i^{(k)}}^{\top}  \right\Vert\\
&=\max_{j}  \sum_{i=1}^n \Vert \E[{x_i}^2\mathbf{x} {\mathbf{x}}^{\top}]  \Vert \left( \sum_{k=1}^N\mathbf{ M}^{(k)}_i {\mathbf{ M}_i^{(k)}}^{\top}  \right)_{jj}.
\end{align*}
To further simplify the right hand side, we note that 
\begin{equation*}
\left(\mathbf{ M}^{(k)}_i {\mathbf{ M}_i^{(k)}}^{\top}  \right)_{jj} = \begin{cases}
\widehat{\gamma}^2_{ii} \delta_i^{(k)} &\text{if } j=i\\
\widehat{\gamma}^2_{ij} \delta_i^{(k)}\delta_j^{(k)} &\text{otherwise}.
\end{cases}
\end{equation*}
Thus for $i=j$, 
\begin{equation*}
\left( \sum_{k=1}^N\mathbf{ M}^{(k)}_i {\mathbf{ M}_i^{(k)}}^{\top}  \right)_{ii} 
= \widehat{\gamma}^2_{ii}\sum_{k=1}^{N}  \delta_i^{(k)}= N \widehat{\gamma}^2_{ii} \widehat{p}_{ii} = N \frac{1}{\widehat{p}_{ii}}.
\end{equation*}
If $i \neq j$ and $(i,j) \in \M$, a similar calculation produces
\begin{equation*}
\left( \sum_{k=1}^N\mathbf{ M}^{(k)}_i {\mathbf{ M}_i^{(k)}}^{\top}  \right)_{ij}=N \frac{1}{\widehat{p}_{ij}}.
\end{equation*}
Finally, when $i \neq j$ and $(i,j) \notin \M$, we have that 
\begin{equation*}
\left( \sum_{k=1}^N\mathbf{ M}^{(k)}_i {\mathbf{ M}_i^{(k)}}^{\top}  \right)_{ij}=\sum_{k=1}^N \widehat{ \gamma}_{ii} \widehat{ \gamma}_{ij}\delta_i^{(k)}\delta_j^{(k)} =0,
\end{equation*}
since $\delta_i^{(k)}\delta_j^{(k)}=0$,  $\forall k$.
We conclude by applying Lemma \ref{lemma_auxiliary_1}, 
\begin{align*}
&\left\Vert \sum_{k=1}^N \mathbb{E}\left[ (\mathbf{M}^{(k)} \odot \mathbf{x}^{(k)} {\mathbf{x}^{(k)}}^{\top} )^2  \mid \mathbf{\Delta}\right] \right\Vert \leq \\ 
&C^2 N \Vert \mathbf{\Sigma}\Vert \max_{j}\sum_{i: (i,j)\in \M} \frac{\Sigma_{ii}}{\widehat{p}_{ij}}.
\end{align*}
\paragraph{Bounding the maximum spectral norm} We follow the same strategy as in the proof of Theorem \ref{th_expected_error} and bound the operator norm by the Frobenius norm. We follow this with Minkowski's inequality leading to the sequence of inequalities
\begin{align*}
&\E \left[  \Vert\mathbf{M}^{(k)} \odot \mathbf{x}^{(k)} {\mathbf{x}^{(k)}}^{\top} \Vert^{2r}  \mid \mathbf{\Delta} \right]^{1/r} \\
&\leq \sum_{i,j} \E\left[\vert x_i^{(k)} x_j^{(k)} \vert^{2r}  \mid \mathbf{\Delta} \right]^{1/r}  \vert M_{ij}^{(k)} \vert^{2}  \\
&= \sum_{i,j} \E\left[\vert x_i^{(k)} x_j^{(k)} \vert^{2r} \right]^{1/r}  \vert M_{ij}^{(k)} \vert^{2}  \\
&\leq \sum_{(i,j) \in \M} \E \left[\vert x_i\vert^{4r}\right]^{1/2r} \E \left[\vert x_j\vert^{4r}\right]^{1/2r} \frac{1}{\widehat{p}^2_{ij}}\\
&\leq (4rc^2)^2  \sum_{(i,j) \in \M} \Vert x_i \Vert_{\Psi_2}^2 \Vert x_j\Vert_{\Psi_2}^2\frac{1}{\widehat{p}^2_{ij}} \\
&\leq  r^2(4c^2/c_0)^2 \sum_{(i,j) \in \M} \frac{\Sigma_{ii}\Sigma_{jj}}{\widehat{p}^2_{ij}}.
\end{align*}
The desired bound is obtained after  setting $r = \log(N)/2$.
%
%
\subsection{Proof of Theorem \ref{th_expected_error_martin} }
We will use an  inequality due to \cite{juditsky2008large}, from where we state the required result in the Lemma below.
\begin{lemma}
	Let $\mathbf{W}_T = \sum_{t=1}^T\tilde{\mathbf{Z}}_t$ be a matrix valued martingale satisfying $\E_{t-1}[\tilde{\mathbf{Z}}_t]=\mathbf{0}$, and $\E[\Vert \tilde{\mathbf{Z}}_t  \Vert^2] \leq \sigma^2_t < \infty$, then
	\begin{equation*}
	\E[\Vert \mathbf{W}_T \Vert^2] \leq c \log(n+1) \sum_{t=1}^T \sigma_t^2.
	\end{equation*}
\end{lemma} 
Taking $\tilde{\mathbf{Z}}_t =\mathbf{Z}_t- N_t\mathbf{\Sigma}$ we just need to bound $\E[\Vert \tilde{\mathbf{Z}}_t  \Vert^2]$. We can do that by conditioning on previous observations and applying  Theorem \ref{th_expected_error}, thus leading to
\begin{align*}
&\E_{t-1}[\Vert \tilde{\mathbf{Z}}_t  \Vert^2] \leq C^2\Vert\mathbf{\Sigma}\Vert^2 N_t^2 \times  \\
&\left[ \sqrt{\frac{8e\log(n) {{\srank}_{\min}}(\mathbf{\Sigma},\mathbf{P}^{(t)}) }{N_t}} + \frac{8e\log(n)\log(N_t){\srank_2}(\mathbf{\Sigma},\mathbf{P}^{(t)})}{N_t}\right]^2 \\
&\leq  2C^2\Vert\mathbf{\Sigma}\Vert^2\times   \\
& \left[8e\log(n){{\srank}_{\min}}(\mathbf{\Sigma},\mathbf{P}^{(t)}) N_t + (8e\log(n)\log(N_t))^2{{\srank}_2}^2(\mathbf{\Sigma},\mathbf{P}^{(t)})\right].
\end{align*}
Then $\sigma_t^2$ is equal to the expectation of the above expression.
We conclude by  adding terms from $t=1$ to $t=T$,   dividing by $(\sum_{t=1}^T N_t)^2$ and computing square root.
\subsection{Proof of Lemma \ref{lemma_auxiliary_1}}
\begin{lemma}\label{lemma_auxiliary_1}
	If  assumption \ref{assum_trace} is true, then 
	\begin{equation*}
	\Vert \E[{x_i}^2\mathbf{x} {\mathbf{x}}^{\top}]  \Vert \leq  C^2 \Sigma_{ii} \Vert \mathbf{\Sigma}\Vert.
	\end{equation*}
\end{lemma}
The proof follows from  the variational definition of operator norm, the Cauchy-Schwartz inequality and Assumption \ref{assum_trace}, leading to the sequence of inequalities below
		\begin{align*}
		\Vert \E[{x_i}^2\mathbf{x} {\mathbf{x}}^{\top}]  \Vert &=\sup_{\Vert \mathbf{u}\Vert_2=1} \mathbf{u}^{\top} \E[{x_i}^2\mathbf{x} {\mathbf{x}}^{\top} ] \mathbf{u} \\
		&= \sup_{\Vert \mathbf{u}\Vert_2=1}  \E[{x_i}^2( {\mathbf{x}}^{\top}\mathbf{u})^2]\\
		&\leq \sup_{\Vert \mathbf{u}\Vert_2=1}  \E[{x_i}^4]^{1/2}\E[( {\mathbf{x}}^{\top}\mathbf{u})^4] ^{1/2}\\
		&\leq 4c^2 \Vert x_i \Vert_{\Psi_2}^2 4c^2  \sup_{\Vert \mathbf{u}\Vert_2=1}  \Vert {\mathbf{x}}^{\top}\mathbf{u} \Vert_{\Psi_2}^2 \\
		&\leq (4c^2/c_0)^2 \Sigma_{ii} \sup_{\Vert \mathbf{u}\Vert_2=1} \mathbf{u}^{\top} \mathbf{\Sigma u} \\
		&=  C^2 \Sigma_{ii} \Vert \mathbf{\Sigma}\Vert.
		\end{align*}
%
\section{Conclusion}
\label{sec:miss_conc}
This paper studied covariance estimation with  various types of  missing observations. We first considered observations missing completely at random (MCAR), where variables are observed according to  non uniform Bernoulli random variables, independent of the data. We studied  unbiased covariance estimators and obtained new  error bounds in operator norm that characterize their performance in terms of the scaled effective rank. These results  improve upon  previous missing data studies,   and are consistent with bounds obtained for the  complete observations case. 
We also introduced the conditionally MCAR missing data mechanism, where the missing data pattern evolves over time and  may depend on previous observations. We proposed an unbiased estimator and  characterized its estimation error  under this setting.

Since in general,  the missing data distribution is unknown, the estimator (\ref{eq_est_cov_y_mu_mcar2}) is suitable for most applications. Our error bounds  indicate that   this estimator achieves the same  convergence rate than others that have access to the missing data distribution. Our numerical experiments support this claim.
%
%
%
%
%
 %
\bibliographystyle{IEEEtran}

\bibliography{references}

\begin{thebibliography}{10}
\providecommand{\url}[1]{#1}
\csname url@samestyle\endcsname
\providecommand{\newblock}{\relax}
\providecommand{\bibinfo}[2]{#2}
\providecommand{\BIBentrySTDinterwordspacing}{\spaceskip=0pt\relax}
\providecommand{\BIBentryALTinterwordstretchfactor}{4}
\providecommand{\BIBentryALTinterwordspacing}{\spaceskip=\fontdimen2\font plus
\BIBentryALTinterwordstretchfactor\fontdimen3\font minus
  \fontdimen4\font\relax}
\providecommand{\BIBforeignlanguage}[2]{{%
\expandafter\ifx\csname l@#1\endcsname\relax
\typeout{** WARNING: IEEEtran.bst: No hyphenation pattern has been}%
\typeout{** loaded for the language `#1'. Using the pattern for}%
\typeout{** the default language instead.}%
\else
\language=\csname l@#1\endcsname
\fi
#2}}
\providecommand{\BIBdecl}{\relax}
\BIBdecl

\bibitem{zhang2010learning}
Y.~Zhang and J.~G. Schneider, ``Learning multiple tasks with a sparse
  matrix-normal penalty,'' in \emph{Advances in Neural Information Processing
  Systems}, 2010, pp. 2550--2558.

\bibitem{schafer2005shrinkage}
J.~Sch{\"a}fer and K.~Strimmer, ``A shrinkage approach to large-scale
  covariance matrix estimation and implications for functional genomics,''
  \emph{Statistical applications in genetics and molecular biology}, vol.~4,
  no.~1, 2005.

\bibitem{bullmore_complex_2009}
E.~Bullmore and O.~Sporns, ``Complex brain networks: graph theoretical analysis
  of structural and functional systems,'' \emph{Nature Reviews Neuroscience},
  vol.~10, no.~3, pp. 186--198, 2009.

\bibitem{bai2011estimating}
J.~Bai and S.~Shi, ``Estimating high dimensional covariance matrices and its
  applications,'' \emph{Annals of Economics and Finance}, vol.~12, no.~2, pp.
  199--215, 2011.

\bibitem{stoica2011spice}
P.~Stoica, P.~Babu, and J.~Li, ``Spice: A sparse covariance-based estimation
  method for array processing,'' \emph{IEEE Transactions on Signal Processing},
  vol.~59, no.~2, pp. 629--638, 2011.

\bibitem{dabov2009bm3d}
K.~Dabov, A.~Foi, V.~Katkovnik, and K.~Egiazarian, ``Bm3d image denoising with
  shape-adaptive principal component analysis,'' in \emph{SPARS'09-Signal
  Processing with Adaptive Sparse Structured Representations}, 2009.

\bibitem{pyatykh2013image}
S.~Pyatykh, J.~Hesser, and L.~Zheng, ``Image noise level estimation by
  principal component analysis,'' \emph{IEEE transactions on image processing},
  vol.~22, no.~2, pp. 687--699, 2013.

\bibitem{chen2010shrinkage}
Y.~Chen, A.~Wiesel, Y.~C. Eldar, and A.~O. Hero, ``Shrinkage algorithms for
  mmse covariance estimation,'' \emph{IEEE Transactions on Signal Processing},
  vol.~58, no.~10, pp. 5016--5029, 2010.

\bibitem{egilmez2017graph}
H.~E. Egilmez, E.~Pavez, and A.~Ortega, ``Graph learning from data under
  laplacian and structural constraints,'' \emph{IEEE Journal of Selected Topics
  in Signal Processing}, vol.~11, no.~6, pp. 825--841, 2017.

\bibitem{lounici2014high}
K.~Lounici, ``High-dimensional covariance matrix estimation with missing
  observations,'' \emph{Bernoulli}, vol.~20, no.~3, pp. 1029--1058, 2014.

\bibitem{bunea2015sample}
F.~Bunea and L.~Xiao, ``On the sample covariance matrix estimator of reduced
  effective rank population matrices, with applications to fpca,''
  \emph{Bernoulli}, vol.~21, no.~2, pp. 1200--1230, 2015.

\bibitem{koltchinskii2017concentration}
V.~Koltchinskii and K.~Lounici, ``Concentration inequalities and moment bounds
  for sample covariance operators,'' \emph{Bernoulli}, vol.~23, no.~1, pp.
  110--133, 2017.

\bibitem{cai2016minimax}
T.~T. Cai and A.~Zhang, ``Minimax rate-optimal estimation of high-dimensional
  covariance matrices with incomplete data,'' \emph{Journal of multivariate
  analysis}, vol. 150, pp. 55--74, 2016.

\bibitem{park2019non}
S.~Park and J.~Lim, ``Non-asymptotic rate for high-dimensional covariance
  estimation with non-independent missing observations,'' \emph{Statistics \&
  Probability Letters}, 2019.

\bibitem{little2014statistical}
R.~J. Little and D.~B. Rubin, \emph{Statistical analysis with missing
  data}.\hskip 1em plus 0.5em minus 0.4em\relax John Wiley \& Sons, 2014, vol.
  333.

\bibitem{kolar2012consistent}
M.~Kolar and E.~P. Xing, ``Consistent covariance selection from data with
  missing values,'' in \emph{Proceedings of the 29th International Conference
  on Machine Learning (ICML)}, 2012, pp. 551--558.

\bibitem{yang2015sparse}
C.~Yang, D.~Robinson, and R.~Vidal, ``Sparse subspace clustering with missing
  entries,'' in \emph{International Conference on Machine Learning}, 2015, pp.
  2463--2472.

\bibitem{hunt2018multi}
\BIBentryALTinterwordspacing
X.~J. Hunt, S.~Emrani, I.~K. Kabul, and J.~Silva, ``Multi-task learning with
  incomplete data for healthcare,'' in \emph{KDD Workshop on Machine Learning
  for Medicine and Healthcare}, 2018. [Online]. Available:
  \url{https://arxiv.org/abs/1807.02442}
\BIBentrySTDinterwordspacing

\bibitem{gonen2016subspace}
A.~Gonen, D.~Rosenbaum, Y.~C. Eldar, and S.~Shalev-Shwartz, ``Subspace learning
  with partial information,'' \emph{The Journal of Machine Learning Research},
  vol.~17, no.~1, pp. 1821--1841, 2016.

\bibitem{chi2013nearest}
Y.~Chi, ``Nearest subspace classification with missing data,'' in
  \emph{Signals, Systems and Computers, 2013 Asilomar Conference on}.\hskip 1em
  plus 0.5em minus 0.4em\relax IEEE, 2013, pp. 1667--1671.

\bibitem{chao2017thesis}
Y.-H. Chao, ``Compression of signal on graphs with the application to image and
  video coding,'' Ph.D. dissertation, University of Southern California, 2017.

\bibitem{loh2012high}
P.-L. Loh and M.~J. Wainwright, ``High dimensional regression with noisy and
  missing data: provable guarantees with non-convexity,'' \emph{The Annals of
  Statistics}, vol.~40, no.~3, pp. 1637--1664, 2012.

\bibitem{jurczak2017spectral}
K.~Jurczak and A.~Rohde, ``Spectral analysis of high-dimensional sample
  covariance matrices with missing observations,'' \emph{Bernoulli}, vol.~23,
  no.~4A, pp. 2466--2532, 2017.

\bibitem{lounici2013sparse}
K.~Lounici, ``Sparse principal component analysis with missing observations,''
  in \emph{High dimensional probability VI}.\hskip 1em plus 0.5em minus
  0.4em\relax Springer, 2013, pp. 327--356.

\bibitem{stadler2012missing}
N.~St{\"a}dler and P.~B{\"u}hlmann, ``Missing values: sparse inverse covariance
  estimation and an extension to sparse regression,'' \emph{Statistics and
  Computing}, vol.~22, no.~1, pp. 219--235, 2012.

\bibitem{veeramachaneni2005active}
S.~Veeramachaneni, F.~Demichelis, E.~Olivetti, and P.~Avesani, ``Active
  sampling for knowledge discovery from biomedical data,'' in \emph{European
  Conference on Principles of Data Mining and Knowledge Discovery}.\hskip 1em
  plus 0.5em minus 0.4em\relax Springer, 2005, pp. 343--354.

\bibitem{lim2005hybrid}
C.-P. Lim, J.-H. Leong, and M.-M. Kuan, ``A hybrid neural network system for
  pattern classification tasks with missing features,'' \emph{IEEE transactions
  on pattern analysis and machine Intelligence}, vol.~27, no.~4, pp. 648--653,
  2005.

\bibitem{asif2016matrix}
M.~T. Asif, N.~Mitrovic, J.~Dauwels, and P.~Jaillet, ``Matrix and tensor based
  methods for missing data estimation in large traffic networks,'' \emph{IEEE
  Transactions on Intelligent Transportation Systems}, vol.~17, no.~7, pp.
  1816--1825, 2016.

\bibitem{deshpande2004model}
A.~Deshpande, C.~Guestrin, S.~R. Madden, J.~M. Hellerstein, and W.~Hong,
  ``Model-driven data acquisition in sensor networks,'' in \emph{Proceedings of
  the Thirtieth international conference on Very large data bases-Volume
  30}.\hskip 1em plus 0.5em minus 0.4em\relax VLDB Endowment, 2004, pp.
  588--599.

\bibitem{gershenfeld2010intelligent}
N.~Gershenfeld, S.~Samouhos, and B.~Nordman, ``Intelligent infrastructure for
  energy efficiency,'' \emph{Science}, vol. 327, no. 5969, pp. 1086--1088,
  2010.

\bibitem{davenport_signal_2010}
M.~A. Davenport, P.~T. Boufounos, M.~B. Wakin, and R.~G. Baraniuk, ``Signal
  {Processing} {With} {Compressive} {Measurements},'' \emph{IEEE Journal of
  Selected Topics in Signal Processing}, vol.~4, no.~2, pp. 445 -- 460, 2010.

\bibitem{dasarathy2015sketching}
G.~Dasarathy, P.~Shah, B.~N. Bhaskar, and R.~D. Nowak, ``Sketching sparse
  matrices, covariances, and graphs via tensor products,'' \emph{IEEE
  Transactions on Information Theory}, vol.~61, no.~3, pp. 1373--1388, 2015.

\bibitem{chen2015exact}
Y.~Chen, Y.~Chi, and A.~J. Goldsmith, ``Exact and stable covariance estimation
  from quadratic sampling via convex programming,'' \emph{IEEE Transactions on
  Information Theory}, vol.~61, no.~7, pp. 4034--4059, 2015.

\bibitem{chen2017toward}
X.~Chen, M.~R. Lyu, and I.~King, ``Toward efficient and accurate covariance
  matrix estimation on compressed data,'' in \emph{International Conference on
  Machine Learning}, 2017, pp. 767--776.

\bibitem{azizyan2018extreme}
M.~Azizyan, A.~Krishnamurthy, and A.~Singh, ``Extreme compressive sampling for
  covariance estimation,'' \emph{IEEE Transactions on Information Theory},
  vol.~64, no.~12, pp. 7613--7635, Dec 2018.

\bibitem{pourkamali2017preconditioned}
F.~Pourkamali-Anaraki and S.~Becker, ``Preconditioned data sparsification for
  big data with applications to pca and k-means,'' \emph{IEEE Transactions on
  Information Theory}, vol.~63, no.~5, pp. 2954--2974, 2017.

\bibitem{pourkamali2016estimation}
F.~Pourkamali-Anaraki, ``Estimation of the sample covariance matrix from
  compressive measurements,'' \emph{IET Signal Processing}, vol.~10, no.~9, pp.
  1089--1095, 2016.

\bibitem{anaraki2014memory}
F.~P. Anaraki and S.~Hughes, ``Memory and computation efficient pca via very
  sparse random projections,'' in \emph{Proceedings of the 31st International
  Conference on Machine Learning (ICML-14)}, 2014, pp. 1341--1349.

\bibitem{naghshvar2013active}
M.~Naghshvar and T.~Javidi, ``Active sequential hypothesis testing,'' \emph{The
  Annals of Statistics}, vol.~41, no.~6, pp. 2703--2738, 2013.

\bibitem{dasarathy2016active}
G.~Dasarathy, A.~Singh, M.-F. Balcan, and J.~H. Park, ``Active learning
  algorithms for graphical model selection,'' in \emph{Artificial Intelligence
  and Statistics (AISTATS)}, 2016, pp. 1356--1364.

\bibitem{scarlett2017lower}
J.~Scarlett and V.~Cevher, ``Lower bounds on active learning for graphical
  model selection,'' in \emph{The 20th International Conference on Artificial
  Intelligence and Statistics (AISTATS)}, 2017.

\bibitem{pavez2018active}
E.~Pavez and A.~Ortega, ``Active covariance estimation by random sub-sampling
  of variables,'' in \emph{2018 IEEE International Conference on Acoustics,
  Speech and Signal Processing (ICASSP)}, April 2018, pp. 4034--4038.

\bibitem{marlin2008missing}
B.~Marlin, ``Missing data problems in machine learning,'' Ph.D. dissertation,
  University of Toronto, 2008.

\bibitem{laird1988missing}
N.~M. Laird, ``Missing data in longitudinal studies,'' \emph{Statistics in
  medicine}, vol.~7, no. 1-2, pp. 305--315, 1988.

\bibitem{zheng2002active}
Z.~Zheng and B.~Padmanabhan, ``On active learning for data acquisition,'' in
  \emph{Data Mining, 2002. ICDM 2003. Proceedings. 2002 IEEE International
  Conference on}.\hskip 1em plus 0.5em minus 0.4em\relax IEEE, 2002, pp.
  562--569.

\bibitem{saar2009active}
M.~Saar-Tsechansky, P.~Melville, and F.~Provost, ``Active feature-value
  acquisition,'' \emph{Management Science}, vol.~55, no.~4, pp. 664--684, 2009.

\bibitem{melville2005economical}
P.~Melville, F.~Provost, M.~Saar-Tsechansky, and R.~Mooney, ``Economical active
  feature-value acquisition through expected utility estimation,'' in
  \emph{Proceedings of the 1st international workshop on Utility-based data
  mining}.\hskip 1em plus 0.5em minus 0.4em\relax ACM, 2005, pp. 10--16.

\bibitem{chakraborty2013active}
S.~Chakraborty, J.~Zhou, V.~Balasubramanian, S.~Panchanathan, I.~Davidson, and
  J.~Ye, ``Active matrix completion,'' in \emph{Data Mining (ICDM), 2013 IEEE
  13th International Conference on}.\hskip 1em plus 0.5em minus 0.4em\relax
  IEEE, 2013, pp. 81--90.

\bibitem{chen2012masked}
R.~Y. Chen, A.~Gittens, and J.~A. Tropp, ``The masked sample covariance
  estimator: an analysis using matrix concentration inequalities,''
  \emph{Information and Inference: A Journal of the IMA}, vol.~1, no.~1, pp.
  2--20, 2012.

\bibitem{vershynin2016high}
R.~Vershynin, \emph{High Dimensional Probability: An Introduction with
  Applications in Data Science}.\hskip 1em plus 0.5em minus 0.4em\relax
  Cambridge University Press, 2018.

\bibitem{durrett2010probability}
R.~Durrett, \emph{Probability: theory and examples}.\hskip 1em plus 0.5em minus
  0.4em\relax Cambridge university press, 2010.

\bibitem{juditsky2008large}
A.~Juditsky and A.~S. Nemirovski, ``Large deviations of vector-valued
  martingales in 2-smooth normed spaces,'' \emph{arXiv preprint
  arXiv:0809.0813}, 2008.

\end{thebibliography}

\appendix
\subsection{Covariance estimation with unknown population mean}  
\label{sec:miss_non_zero_mean}
In this section we derive unbiased estimators when  the population mean $\boldsymbol\mu$ is {\em unknown}, and the missing data is MCAR. We consider  two scenarios:
\begin{enumerate}
	\item { Known missing data distribution}: the estimators depend on the observed data $\mathbf{y}^{(k)}$ for $k \in [N]$, the population observation probabilities $\mathbf{p}$, $\mathbf{P}$, and their entrywise inverses $\boldsymbol{\gamma}$ and $\mathbf{\Gamma}$ respectively.
	\item { Unknown missing data distribution}: the estimators depend on the observed data $\mathbf{y}^{(k)}$ for $k \in [N]$, the sample estimators for the observation probabilities $\hat{\mathbf{p}}$, $\widehat{\mathbf{P}}$, and their entrywise inverses $\hat{\boldsymbol{\gamma}}$ and $\widehat{\mathbf{\Gamma}}$ respectively.
\end{enumerate}
Some quantities related to the missing data mechanism are given next. The sample observation probabilities 
\begin{equation}
	\hat{\mathbf{p}} = \frac{1}{N} \sum_{k=1}^{N}\boldsymbol\delta^{(k)}.
\end{equation}
The sample observation  probability matrix is 
\begin{equation}
	\widehat{\mathbf{P}} = \frac{1}{N} \sum_{k=1}^{N}\boldsymbol\delta^{(k)} {\boldsymbol\delta^{(k)}}^{\top}.
\end{equation}
When $\boldsymbol\delta^{(k)}$ has dependent Bernoulli entries, we have that
\begin{align}
	\E[\widehat{\mathbf{P}}] &= \mathbf{P}, \\
	\E[\hat{\mathbf{p}}] &= \mathbf{p} = \diag(\mathbf{P}).
\end{align}
Recall that $\mathbf{\Gamma}$ and $\widehat{\mathbf{\Gamma}}$ are the entry-wise inverses of $\mathbf{P}$ and $\widehat{\mathbf{P}}$ respectively. The  entry-wise inverses of $\mathbf{p}$ and $\widehat{\mathbf{p}}$ are defined as $\boldsymbol\gamma$ and $\hat{\boldsymbol\gamma}$ respectively. 

It will also be useful to have estimates for the mean and covariance of the incomplete observations. The sample mean of $\mathbf{y}$ is
\begin{equation}
	\bar{\mathbf{y}} = \frac{1}{N} \sum_{k=1}^{N} \mathbf{y}^{(k)}.
\end{equation}
This estimator has expectation given by
\begin{equation}
\E[\bar{\mathbf{y}} \mid \mathbf{\Delta}] = \hat{\mathbf{p}} \odot \boldsymbol\mu, \textnormal{ and } \quad \E[\bar{\mathbf{y}} ] = \mathbf{p} \odot \boldsymbol\mu.
\end{equation}
The covariance of $\mathbf{y}$ can be estimated with
\begin{equation}
 \Cy = \frac{1}{N-1}  \sum_{k=1}^N (\mathbf{ y}^{(k)} -  \bar{\mathbf{y}} )   (\mathbf{ y}^{(k)} -  \bar{\mathbf{y}} )^{\top}.
\end{equation}
When $\mathbf{P}$ is known, $\Cy$ is unbiased for $\Cov(\mathbf{y})$, thus
\begin{equation}\label{eq_sample_covarianze_nonzero_mean_unknown}
	\E[\Cy] = \mathbf{ \Sigma } \odot \mathbf{P} + \left(\mathbf{P} - \mathbf{p}\mathbf{p}^{\top}\right) \odot \boldsymbol\mu {\boldsymbol\mu}^{\top}=\Cov(\mathbf{y}).
\end{equation}
Note that since the goal is to estimate $\mathbf{ \Sigma }$, we could attempt to correct  $\Cy$ by removing  the effect of $\mathbf{P}$ and $\boldsymbol\mu$. A simpler strategy is to   build estimators using $\bar{\mathbf{y}}$, and
\begin{equation}
\Ry  = \frac{1}{N}\sum_{k=1}^N \mathbf{ y}^{(k)}   {\mathbf{ y}^{(k)}}^{\top}.
\end{equation}
Before deriving the covariance estimators,  we introduce some useful properties of $\bar{\mathbf{y}}$ and $\Ry$. Proofs are provided in the subsequent appendix.
%
%
%
\begin{lemma}\label{lemma_Ry}
The quantity $\Ry$ obeys
\begin{align}
	\E\left[\Ry  \mid \mathbf{\Delta } \right] &=  \widehat{\mathbf{P}} \odot (\mathbf{\Sigma} + \boldsymbol\mu \boldsymbol\mu^{\top})\\
		\E\left[\Ry   \right]&= {\mathbf{P}} \odot (\mathbf{\Sigma} + \boldsymbol\mu \boldsymbol\mu^{\top}).
\end{align}
\end{lemma}
\begin{lemma}\label{lemma_y_P1}
When $\mathbf{P}$ is known we have  
\begin{equation*}
\E[\bar{\mathbf{y}} \bar{\mathbf{y}}^{\top}]=\frac{1}{N} \left( \mathbf{ \Sigma } \odot \mathbf{P} + \left(\mathbf{P} - \mathbf{p}\mathbf{p}^{\top}\right) \odot \boldsymbol\mu {\boldsymbol\mu}^{\top} \right)  +  \mathbf{p}\mathbf{p}^{\top} \odot \boldsymbol\mu {\boldsymbol\mu}^{\top}.
\end{equation*}
\end{lemma}
\begin{lemma}\label{lemma_y_P2}
When $\mathbf{P}$ is unknown we have  
\begin{equation*}
\E[\bar{\mathbf{y}} \bar{\mathbf{y}}^{\top} \mid \mathbf{\Delta}] = \hat{\mathbf{p}} \hat{\mathbf{p}}^{\top } \odot \boldsymbol\mu \boldsymbol\mu^{\top}  + \frac{1}{N} \widehat{\mathbf{P}} \odot \mathbf{\Sigma}.
\end{equation*}
\end{lemma}
\subsubsection{Covariance estimation with known $\mathbf{P}$}
According to   Lemma \ref{lemma_Ry}, $\Ry \odot \mathbf{\Gamma}$ 
is an unbiased estimator for $\mathbf{\Sigma} + \boldsymbol\mu \boldsymbol\mu^{\top}$.  Our strategy is thus to find an unbiased estimator for $ \boldsymbol\mu \boldsymbol\mu^{\top}$ and subtract it from $\Ry \odot \mathbf{\Gamma}$.  
 Since Lemma  \ref{lemma_Ry} and \ref{lemma_y_P1}  imply
\begin{equation}
		\E\left[N^2 \bar{\mathbf{y}} \bar{\mathbf{y}}^{\top} - N\Ry   \right] = N(N-1) \mathbf{p}\mathbf{p}^{\top} \odot \boldsymbol\mu {\boldsymbol\mu}^{\top}, 
\end{equation}
an unbiased estimator for $\boldsymbol\mu {\boldsymbol\mu}^{\top}$ is given by
\begin{equation}
	\frac{1}{N-1}\left( N \bar{\mathbf{y}} \bar{\mathbf{y}}^{\top} - \Ry \right)  \odot \boldsymbol\gamma \boldsymbol\gamma^{\top}.
\end{equation}
Therefore, an  unbiased estimator for $\mathbf{\Sigma}$ is 
\begin{equation}
	\Ry \odot \mathbf{\Gamma} - \frac{1}{N-1}\left( N \bar{\mathbf{y}} \bar{\mathbf{y}}^{\top} - \Ry \right)  \odot \boldsymbol\gamma \boldsymbol\gamma^{\top}
\end{equation}
\subsubsection{Covariance estimation with unknown $\mathbf{P}$}
We can derive an  estimator following a similar procedure, however here we estimate $\boldsymbol\mu \boldsymbol\mu^{\top}$ differently. The proposed covariance estimator corresponds to
\begin{equation}
		\Ry \odot \widehat{\mathbf{\Gamma}} - ( N \bar{\mathbf{y}} \bar{\mathbf{y}}^{\top} - \Ry ) \odot\mathbf{ \Psi}.
\end{equation}
Where $\mathbf{ \Psi}$ is the Hadamard inverse of  $\mathbf{\Theta}=N\hat{\mathbf{p}} \hat{\mathbf{p}}^{\top } - \hat{\mathbf{P}}$. Unbiasedness can be verified by applying  Lemma \ref{lemma_Ry} and \ref{lemma_y_P2}. A necessary condition for $\mathbf{\Theta}$ to be Hadamard invertible is that $\sum_{k =1}^N\delta_i^{(k)} \geq 2$ for all $i \in [n]$. This result is proven in the next appendix.
\subsubsection{Comparison with complete observations}
Note that when there is no missing data, both covariance estimators from this section become
\begin{equation}
\frac{1}{N-1}\sum_{k=1}^N (\mathbf{y}^{(k)} - \bar{\mathbf{y}})(\mathbf{y}^{(k)} - \bar{\mathbf{y}})^{\top}.
\end{equation}
\subsection{Additional proofs}
\subsubsection{Proof of Lemma \ref{lemma_y_P1} }
We may use $\boldsymbol\mu_{\mathbf{y}} = \mathbf{p} \odot \boldsymbol\mu $, then
\begin{align*}
&N^2	\bar{\mathbf{y}} \bar{\mathbf{y}}^{\top}= \sum_{k=1}^N \sum_{l=1}^N\mathbf{y}^{(k)} {\mathbf{y}^{(l)}}^{\top} \\
&= \sum_{k=1}^N \sum_{l=1}^N (\mathbf{y}^{(k)}- \boldsymbol\mu_{\mathbf{y}}) (\mathbf{y}^{(l)}- \boldsymbol\mu_{\mathbf{y}})^{\top} + N^2  \boldsymbol\mu_{\mathbf{y}} \boldsymbol\mu_{\mathbf{y}}^{\top}\\
&+ \sum_{k=1}^N \sum_{l=1}^N [(\mathbf{y}^{(k)}- \boldsymbol\mu_{\mathbf{y}}) { \boldsymbol\mu_{\mathbf{y}}}^{\top} +  \boldsymbol\mu_{\mathbf{y}} (\mathbf{y}^{(k)}- \boldsymbol\mu_{\mathbf{y}})^{\top}]  \\
&= N \Sy  + 2\sum_{k>l}^N (\mathbf{y}^{(k)}- \boldsymbol\mu_{\mathbf{y}}) (\mathbf{y}^{(l)}- \boldsymbol\mu_{\mathbf{y}})^{\top} + N^2  \boldsymbol\mu_{\mathbf{y}} \boldsymbol\mu_{\mathbf{y}}^{\top} \\
&+ N(	\bar{\mathbf{y}} -\boldsymbol\mu_{\mathbf{y}} ) \boldsymbol\mu_{\mathbf{y}}^{\top } + N\boldsymbol\mu_{\mathbf{y}}(	\bar{\mathbf{y}} -\boldsymbol\mu_{\mathbf{y}} )^{\top}.
\end{align*}
We conclude by applying expectation and using the fact that $\E[\Sy] = \Cov(\mathbf{y})$.
\subsubsection{Proof of Lemma \ref{lemma_y_P2}}
 We have the identity
\begin{align*}
N^2	\bar{\mathbf{y}} \bar{\mathbf{y}}^{\top}&= \sum_{k=1}^N \sum_{l=1}^N\mathbf{y}^{(k)} {\mathbf{y}^{(l)}}^{\top}\\
&=  \sum_{k \neq l}\mathbf{y}^{(k)} {\mathbf{y}^{(l)}}^{\top} +  \sum_{k=1}^N \mathbf{y}^{(k)} {\mathbf{y}^{(k)}}^{\top}.
\end{align*}
Applying conditional expectation leads to
\begin{align*}
	&\E[N^2 \bar{\mathbf{y}} \bar{\mathbf{y}}^{\top} \mid \mathbf{\Delta}] = \sum_{k \neq l}\boldsymbol\delta^{(k)} {\boldsymbol\delta^{(l)}}^{\top} \odot \boldsymbol\mu \boldsymbol\mu^{\top} \\
	&+  \sum_{k=1}^N \boldsymbol\delta^{(k)} {\boldsymbol\delta^{(k)}}^{\top} \odot (\mathbf{\Sigma} + \boldsymbol\mu \boldsymbol\mu^{\top}) \\
	&=\sum_{k=1}^N \sum_{l=1}^N\boldsymbol\delta^{(k)} {\boldsymbol\delta^{(l)}}^{\top} \odot \boldsymbol\mu \boldsymbol\mu^{\top} - \sum_{k=1}^N \boldsymbol\delta^{(k)} {\boldsymbol\delta^{(k)}}^{\top} \odot \boldsymbol\mu \boldsymbol\mu^{\top}\\
	&+  \sum_{k=1}^N \boldsymbol\delta^{(k)} {\boldsymbol\delta^{(k)}}^{\top} \odot (\mathbf{\Sigma} + \boldsymbol\mu \boldsymbol\mu^{\top}) \\
	&=N^2 \hat{\mathbf{p}} \hat{\mathbf{p}}^{\top } \odot \boldsymbol\mu \boldsymbol\mu^{\top}  + N \widehat{\mathbf{P}} \odot \mathbf{\Sigma}.
\end{align*}
\subsubsection{Hadamard invertibility of $\mathbf{\Theta}$}
The $ij$ entry of $\mathbf{\Theta}$ obeys
\begin{equation}\label{eq_entry_correction_term_final}
N\theta_{ij}=  \sum_{k =1}^N\delta_i^{(k)} \sum_{k =1}^N\delta_j^{(k)} - \sum_{k =1}^N\delta_i^{(k)}\delta_j^{(k)} \geq 0.
\end{equation}
When $i = j$, a necessary and sufficient condition for invertibility is that  for all $i \in [n]$,  $ \sum_{k =1}^N\delta_i^{(k)} \geq 1$.
For $i \neq j$,  $\theta_{ij}>0$ if any of the following conditions holds
\begin{enumerate}
\item  $\sum_{k =1}^N\delta_i^{(k)} = \sum_{k =1}^N\delta_j^{(k)} = 1$, and $\delta_i^{(k)} \delta_j^{(k)} = 0$ for all $k$, 
\item $\sum_{k =1}^N\delta_i^{(k)} \geq 1$ and $\sum_{k =1}^N\delta_j^{(k)} \geq 2$, 
\item $\sum_{k =1}^N\delta_i^{(k)} \geq 2$ and $\sum_{k =1}^N\delta_j^{(k)} \geq 1$.
\end{enumerate}
Case 1) can be verified directly. Case 2) and 3) can be proven using the same argument. Let us consider case 2) first.  Since $0 \leq \delta_j^{(k)} \leq 1$ we have that
\begin{equation}
 \sum_{k =1}^N\delta_i^{(k)}\delta_j^{(k)} \leq  \sum_{k =1}^N\delta_i^{(k)},
\end{equation}
which combined with  $ \sum_{k =1}^N\delta_j^{(k)} >1$ implies that $\theta_{ij}>0$. 

 %
\end{document}